\documentclass[a4paper]{amsart}

\usepackage[T1]{fontenc}
\usepackage[latin1]{inputenc}
\usepackage{amssymb}
\usepackage{euscript}
\usepackage{amsxtra}
\usepackage{ae}
\usepackage[all]{xy}
\usepackage{enumerate}
\usepackage[lite]{amsrefs}
\usepackage{mathrsfs}
\usepackage[german,british]{babel}

\theoremstyle{plain}
\newtheorem{thm}{Theorem}
\newtheorem{lemma}[thm]{Lemma}
\newtheorem{prop}[thm]{Proposition}
\newtheorem{cor}[thm]{Corollary}
\theoremstyle{definition}
\newtheorem{defn}[thm]{Definition}
\theoremstyle{remark}
\newtheorem{rmk}[thm]{Remark}

\newcommand*{\C}{\mathbb C}
\newcommand*{\Z}{\mathbb Z}
\newcommand*{\R}{\mathbb R}
\newcommand*{\N}{\mathbb N}

\newcommand*{\vvr}{\mathfrak{\bar{c}}^{\mathfrak{red}}}
\newcommand*{\bvr}{\mathfrak{c}^{\mathfrak{red}}}

\newcommand*{\bbbarred}{\bar{\mathfrak{B}}^{\mathfrak{red}}}
\newcommand*{\bbred}{\mathfrak{B}^{\mathfrak{red}}}
\newcommand*{\bbbar}{\bar{\mathfrak{B}}}
\newcommand*{\ADir}{{\mathsf{P}}}
\newcommand*{\Dirac}{\mathsf{D}}

\newcommand*{\KK}{\mathrm{KK}}
\newcommand*{\RKK}{\mathrm{RKK}}
\newcommand*{\K}{\mathrm{K}}
\newcommand*{\Ktop}{\mathrm{K}^{\mathrm{top}}}
\newcommand*{\KX}{\mathrm{KX}}

\newcommand*{\fin}{{\mathrm{fin}}}
\newcommand*{\Jet}{{\mathrm{J}}}
\newcommand*{\Spin}{{\mathrm{Spin}}}
\newcommand*{\ID}{{\mathrm{id}}}

\newcommand*{\GEN}{\mathcal{G}}
\newcommand*{\Tri}{\mathcal{T}}

\DeclareMathOperator{\Res}{Res}
\DeclareMathOperator{\Ind}{Ind}

\DeclareMathOperator{\supp}{supp}
\DeclareMathOperator{\Diff}{Diff}

\newcommand*{\Hilm}{{\mathcal{E}}}
\newcommand*{\CI}{{\mathcal{CI}}}
\newcommand*{\EG}{{\mathcal{E}G}}

\newcommand*{\Rips}{\mathscr{P}}
\newcommand*{\Subrips}{\mathrm{P}}
\newcommand*{\X}{\mathscr{X}}
\newcommand*{\Y}{\mathscr{Y}}

\newcommand*{\coa}{\mathscr{C}}
\newcommand*{\cocoa}{\mathscr{CC}}

\newcommand*{\Mult}{{\mathcal{M}}}
\newcommand*{\Cmax}{C^*_{\mathrm{max}}}
\newcommand*{\Cred}{C^*_{\mathrm{r}}}

\newcommand*{\Comp}{{\mathbb{K}}}
\newcommand*{\Bound}{{\mathbb{B}}}
\newcommand*{\Calk}{\mathcal{Q}}
\newcommand*{\brd}{-\hspace{0pt}}
\newcommand*{\nbd}{\nobreakdash-\hspace{0pt}}

\newcommand*{\abs}[1]{\lvert#1\rvert}
\newcommand*{\ket}[1]{\lvert#1\rangle}
\newcommand*{\bra}[1]{\langle#1\rvert}
\newcommand*{\norm}[1]{\lVert#1\rVert}
\newcommand*{\gen}[1]{\langle#1\rangle}

\newcommand*{\rcross}{\mathbin{\rtimes_{\mathrm{r}}}}
\newcommand*{\cross}{\mathbin{\rtimes}}
\newcommand*{\defeq}{\mathrel{:=}}
\newcommand*{\into}{\rightarrowtail}
\newcommand*{\prto}{\twoheadrightarrow}

\begin{document}

\title[A descent principle for the Dirac dual Dirac method]
  {A descent principle for the\\ Dirac dual Dirac method}

\author{Heath Emerson}
\email{hemerson@math.uni-muenster.de}

\author{Ralf Meyer}
\email{rameyer@math.uni-muenster.de}

\address{\selectlanguage{german}Mathematisches Institut\\
         Westf.\ Wilhelms-Universität Münster\\
         Einsteinstraße 62\\
         48149 Münster\\
         Deutschland}

\begin{abstract}
  Let~$G$ be a torsion free discrete group with a finite dimensional
  classifying space $BG$.  We show that~$G$ has a dual Dirac morphism
  if and only if a certain coarse (co)-assembly map is an isomorphism.
  Hence the existence of a dual Dirac morphism for such~$G$ is a
  metric, that is, coarse, invariant of~$G$.  We get similar results
  for groups with torsion as well.  The framework that we develop is
  also suitable for studying the Lipschitz and proper Lipschitz
  cohomology of Connes, Gromov and Moscovici.
\end{abstract}


\thanks{This research was supported by the EU-Network \emph{Quantum
  Spaces and Noncommutative Geometry} (Contract HPRN-CT-2002-00280)
  and the \emph{Deutsche Forschungsgemeinschaft} (SFB 478).}

\maketitle

\section{Introduction}
\label{sec:intro}

Let~$G$ be a discrete group with finite classifying space $BG$.  The
\emph{Descent Principle} (see~\cite{HigRoe}) asserts that the
Baum-Connes assembly map $\K_*(BG)\to\K_*(\Cred G)$ is injective if
the the coarse Baum-Connes assembly map $\KX_*(G)\to
\K_*\bigl(C^*(\abs{G})\bigr)$ is an isomorphism.  The latter assertion
only involves the large scale geometry of~$G$.  Injectivity of the
Baum-Connes assembly map implies the homotopy invariance of higher
signatures for~$G$.  It has been expected for some time that there
should be a similar descent principle underlying the Dirac dual Dirac
method.  The primary goal of this article is to specify exactly such a
principle.

Let~$G$ be a locally compact group.  There is a canonical coarse
structure on~$G$ that is invariant under right multiplication.  We
write~$\abs{G}$ for~$G$ equipped with this coarse structure.  Let
$\EG$ be a universal proper $G$\nbd{}space.  Consider the map
\begin{equation}  \label{eq:analytic_coassembly}
  p_\EG^*\colon
  \KK^G_*\bigl(\C,C_0(G)\bigr) \to \RKK^G_*\bigl(\EG;\C,C_0(G)\bigr)
\end{equation}
which is induced by the constant map from $\EG$ to a point.  This
article grew out of the observation
that~\eqref{eq:analytic_coassembly} is an invariant of the coarse
space~$\abs{G}$.  More precisely, given any coarse space~$X$, there is
a certain $C^*$\nbd{}algebra $\bvr(X)$ called the \emph{reduced stable
  Higson corona} of~$X$, a certain graded Abelian group $\KX^*(X)$
called the \emph{coarse $\K$\nbd{}theory} of~$X$, and a \emph{coarse
  co-assembly map}
\begin{equation}  \label{coarsecoassemblymapfirstmentionofit}
  \mu_X^*\colon \K_{*+1}\bigl(\bvr(X)\bigr) \to \KX^*(X),
\end{equation}
which is equivalent to~\eqref{eq:analytic_coassembly} for $X=\abs{G}$.
``\emph{Equivalent}'' means that there are isomorphisms between the
sources and targets of $\mu_{\abs{G}}^*$ and $p_\EG^*$ that identify
the two maps in the obvious way.

This article is an expanded and completely rewritten version of the
eprint~\cite{EmersonMeyer1}.  In the meantime, we have published
details concerning the stable Higson corona and the coarse
co\brd{}assembly map in~\cite{EmersonMeyer}.  Thus we shall only
recall these constructions rather briefly here.

The map~\eqref{eq:analytic_coassembly} is closely related to the Dirac
dual Dirac method for proving injectivity of the Baum-Connes assembly
map.  Traditionally, this method requires the existence of a proper
$G$-$C^*$-algebra~$A$ and $d\in\KK^G(A,\C)$, $\eta\in\KK^G(\C,A)$ such
that $\eta\circ d=1_A$ or at least $p_\EG^*(d\circ\eta)=1_\C$
(see~\cite{KasparovSkandalis2}).  In~\cite{MeyerNest} it is shown that
there is a canonical choice for the \emph{Dirac morphism}
$\Dirac\in\KK^G(\ADir,\C)$. It is constructed using general results
on triangulated categories and can be interpreted as a projective
resolution of~$\C$ in the category $\KK^G$ with respect to a certain
localising subcategory.  The Dirac morphism is unique up to
$\KK^G$\brd{}equivalence.  A \emph{dual Dirac morphism} is defined
in~\cite{MeyerNest} as an element $\eta\in\KK^G(\C,\ADir)$ with
$\eta\circ\Dirac=1_\ADir$.  The dual Dirac morphism is unique if it
exists.  If the Dirac dual Dirac method applies to~$G$ in the
(traditional) sense of~\cite{KasparovSkandalis2}, then there is a dual
Dirac morphism in the sense of~\cite{MeyerNest}.  The converse also
holds for many groups~$G$.  For instance, it holds if~$G$ is discrete
and has a finite dimensional model for $\EG$.

It is further shown in~\cite{MeyerNest} that the canonical maps
\begin{equation}  \label{eq:Ktop_with_Dirac}
  \Ktop_*(G,A)
  \overset{\cong}\leftarrow \Ktop_*(G,\ADir \otimes A)
  \overset{\cong}\to \K_*\bigl((\ADir\otimes A) \cross G\bigr)
  \overset{\cong}\to \K_*\bigl((\ADir\otimes A) \rcross G\bigr)
\end{equation}
are all isomorphisms.  Here $\rcross$ and~$\cross$ denote the reduced
and full crossed products, respectively.  Therefore, we may identify
$\Ktop_*(G,A)$ with $\K_*\bigl((\ADir\otimes A)\cross G\bigr)$ and the
Baum-Connes assembly map with the composition
$$
\K_*\bigl((\ADir\otimes A) \cross G\bigr)
\overset{\Dirac_*}\to \K_*(A\cross G)
\to \K_*(A\rcross G).
$$
The isomorphism appearing in equation~\eqref{intro:fundamentalisomorphism} below, 
and many other of our
results, depend on this definition of the Baum-Connes assembly map, and
cannot apparently be proved with the traditional one.

Let~$G$ be a torsion free discrete group with finite classifying space
$BG$.  Our first objective is to prove that such a group~$G$ has a
dual Dirac morphism if and only if~\eqref{eq:analytic_coassembly} is
an isomorphism.  As a result, the Dirac dual Dirac method applies
to~$G$ if and only if the coarse co\brd{}assembly map
$\mu_{\abs{G}}^*$ is an isomorphism.  This is our new Descent
Principle.  We obtain a similar criterion if~$G$ is a torsion free
discrete group with finite \emph{dimensional} $BG$.  This
generalisation requires the coarse co\brd{}assembly map with
coefficients.

As a result, for torsion free discrete groups with finite dimensional
classifying space, the existence of a dual Dirac morphism is a coarse
invariant: if two such groups $G$ and~$G'$ are coarsely equivalent,
then~$G$ has a dual Dirac morphism if and only if~$G'$ has one.

We also provide results for locally compact groups with a
$G$\nbd{}compact model for $\EG$.  For such groups a dual Dirac
morphism exists if and only if the maps
\begin{equation}  \label{eq:analytic_coassembly_H}
  p_\EG^*\colon \KK^G_*\bigl(\C,C_0(G/H)\bigr)
  \to\RKK^G_*\bigl(\EG;\C,C_0(G/H)\bigr)
\end{equation}
are isomorphisms for all compact subgroups $H\subseteq G$.  We can
equip $\bvr(\abs{G})$ with a canonical action of~$G$.  This gives rise
to an \emph{$H$\nbd{}equivariant coarse co\brd{}assembly map}
$$
\mu_{\abs{G},H}^*\colon
\K_{*+1}\bigl(\bvr_H(\abs{G})^H\bigr) \to \KX_H^*(\abs{G}),
$$
for each compact subgroup $H\subseteq G$, which is equivalent
to~\eqref{eq:analytic_coassembly_H}.

From the perspective of Novikov's original conjecture, the Dirac dual
Dirac method is a method of confirming
 homotopy invariance of all
higher signatures for a group~$G$ at once.  Our Descent Principle
shows that the success or failure of this method depends only on the
large scale geometry of~$G$.  Instead, one may attempt to 
establish the homotopy invariance of a single higher signature,
should it fortuitously arise from a particularly geometric construction. 
The stable Higson corona construction provides a good framework for
pursuing this idea, which goes back to Connes, Gromov and Moscovici in
~\cite{CGM2}. We shall use
it to simplify the geometric part of the proof of homotopy invariance
of Gelfand-Fuchs classes, a result of~\cite{CGM2}. More
substantially, 
we introduce a generalization of the notion of \emph{proper Lipschitz
class}, which we term \emph{boundary class}, and show that 
such classes have a number of pleasant properties analogous
to those enjoyed by Lipschitz classes, including 
homotopy invariance.

For the purposes of this discussion, let us agree that a higher
signature for a discrete group~$G$ is a linear map $\tau\colon
\Ktop_*(G)\to\R$.  If $BG$ is a compact space, then a higher
signature~$\tau$ in our sense may be used to assign a real number to
each continuous map $f\colon M\to BG$, where~$M$ is a smooth oriented
manifold, by the formula $f\mapsto \tau(f_*[D^{\mathrm{sig}}_M])$.
Homotopy invariance means the real numbers associated to $f\colon M
\to BG$ and $f\circ\phi\colon N\to BG$ are the same, where~$\phi$ is a
homotopy equivalence $N\to M$ (see~\cite{Ros}).  It is well-known that
a higher signature is homotopy invariant if it factorises through the
analytic assembly map $\mu_*\colon \Ktop_*(G)\to\K_*(\Cmax G)$
(see~\cite{KM}).

If~$A$ is any $G$-$C^*$-algebra, then $\Ktop_*(G,A)$ is a graded
module over the graded commutative ring $\RKK^G_*(\EG;\C,\C)$.  In
particular, this is true of $\Ktop_*(G)$.  Hence any class
$b\in\RKK^G_*(\EG;\C,\C)$ yields a $\K$\nbd{}theoretic higher
signature by the following recipe.  We set $\tau_b(a) =
\mathrm{ind}_G(a\cdot b)$, where $\mathrm{ind}_G\colon
\Ktop_*(G)\to\Z$ is the index map induced by the trivial
representation of~$G$.  If~$b$ lies in the range of the map
$p_\EG^*\colon \KK^G_*(\C,\C)\to\RKK^G_*(\EG;\C,\C)$ then~$\tau_b$ is
homotopy invariant because it factorises through~$\mu_*$.  If~$G$ has
a dual Dirac morphism, then $p_\EG^*$ is split surjective. This yields 
homotopy invariance of all higher signatures at once. 

%By contrast, Connes, Gromov and Moscovici in~\cite{CGM2} describe
%a geometric condition on a \emph{single} higher signature
%that implies its homotopy invariance. Their condition is that 
%the class arise from a certain construction involving proper Lipschitz
%maps from $G$-spaces to Euclidean space. These higher signatures 
%are called Lipschitz cohomology classes. 

By contrast, and following the spirit of 
\cite{CGM2}, we show how to construct individual 
homotopy invariant higher signatures
directly from the stable Higson corona of $G$. The latter is 
in a natural way a
$G$\nbd{}$C^*$\brd{}algebra, so that we may consider the graded group
$\Ktop_*\bigl(G,\bvr(\abs{G})\bigr)$.  Since we use $\Ktop$, the group
structure of~$G$ does not add any additional \emph{analytical}
complications.  There is a canonical map
\begin{equation}  \label{intro:canonicalmap}
\Ktop_{*+1}\bigl(G,\bvr(\abs{G})\bigr)\to \RKK^G_*(\EG;\C,\C).
\end{equation}
We call elements in the range of this map \emph{boundary classes}
for the obvious reason: they depend only on the action of $G$ on its
stable Higson corona. Suppose $G$ is a discrete group with a $G$-compact model for $\EG$. 
Let $\ADir$ be the source of the Dirac morphism. Then there is
 an isomorphism 

\begin{equation}  \label{intro:fundamentalisomorphism}
  \Ktop_{*+1}\bigl(G,\bvr(\abs{G})\bigr) \cong \KK^G_*(\C,\ADir).
\end{equation}

%Boundary classes do indeed give homotopy invariant higher signatures. The source and 
%target of~\eqref{intro:canonicalmap} are modules over the ring
%$\RKK^G(\EG ; \C , \C)$ and the map~\eqref{intro:canonicalmap} is a module
%homomorphism. Consequently the boundary classes form an ideal in 
%the ring $\RKK^G_*(\EG ; \C ,\C)$. All these assertions follow from 
%the the following fundamental result. 

%It is worth noting that this suggests on 
%purely algebraic grounds that boundary classes are rather special 
% even amongst
%those classes producing homotopy invariant signatures.
% For the latter \emph{a priori} certainly
%do not form an ideal. The unit
%$1_\EG\in\RKK^G_*(\EG;\C,\C)$ always yields a homotopy invariant
%higher signature, namely the ordinary signature, assigning to a map
%$f\colon M \to BG$ the Hirzebruch signature of~$M$. By constrast,
%if $1_\EG$ is a boundary class, then all classes are 
%boundary classes. As we show, this actually implies that
% $G$ has a dual Dirac class. 

This result is fundamental for the study of 
boundary classes. As well, it allows us to significantly 
refine our 
Descent Principle for the Dirac dual Dirac method.

Firstly,~\eqref{intro:fundamentalisomorphism}
 implies that any boundary class yields a
 homotopy invariant 
 higher signature. In fact, any boundary class must lie in the 
range of the map $p_\EG^*$. 

Secondly, it implies that 
\emph{every proper Lipschitz class} in the sense of Connes, Gromov
and Moscovici is a boundary class. Therefore the notion of 
boundary class generalizes that of proper Lipschitz class. 

Thirdly, it shows that a dual 
Dirac morphism may be described directly in terms of the topological 
$\K$-theory of the stable Higson corona, namely as an element 
of $\Ktop_1\bigl(G ; \bvr(G)\bigr)$ whose image under the 
map~\eqref{intro:canonicalmap} is $1_\EG$. Actually, 
using our Descent Principle, we conclude that~\eqref{intro:canonicalmap}
is an isomorphism exactly when it is surjective, and this 
occurs exactly when $G$ has a dual Dirac morphism. 

In particular, if $G$ admits a uniform embedding in Hilbert space, 
every class is a boundary class, and~\eqref{intro:canonicalmap}
is an isomorphism. For such $G$ has a dual Dirac morphism by 
\cite{EmersonMeyer}.

The isomorphism~\eqref{intro:fundamentalisomorphism} also
 enables us to 
construct dual Dirac morphisms directly 
from contractible, admissible,
$G$\nbd{}equivariant compactifications $\overline{\EG} =
\EG\cup\partial G$.  This strengthens a result of Nigel Higson
(\cite{Hig2}).  For instance, if~$G$ is Gromov hyperbolic, we may
compactify $\EG$ using the Gromov boundary of~$G$. In particular, 
a dual Dirac morphism can be constructed which, in the 
appropriate sense, extends continuously
to the Gromov boundary of $G$.

The source and 
target of~\eqref{intro:canonicalmap} are modules over the ring
$\RKK^G(\EG ; \C , \C)$. Moreover, 
the map~\eqref{intro:canonicalmap} is a module
homomorphism. Hence
 the boundary classes form an \emph{ideal} in 
the ring $\RKK^G_*(\EG ; \C ,\C)$. This shows that 
boundary classes are rather special,
 even amongst classes in the range of $p_\EG^*$, or even 
more
generally, amongst those classes which yield homotopy invariant
 higher signatures. For the unit 
$1_\EG\in\RKK^G_*(\EG;\C,\C)$ always yields a homotopy invariant
higher signature, namely the ordinary signature, assigning to a map
$f\colon M \to BG$ the Hirzebruch signature of~$M$. On the other hand,
if $1_\EG$ is a boundary class, then, since the 
boundary classes constitute an ideal, all classes are 
boundary classes. As remarked above, this is the case if and only if 
 $G$ has a dual Dirac class.

Finally, since the source of the Dirac morphism and consequently also
the group $\KK^G(\C,\ADir)$ are already built into the definition of
the Baum-Connes assembly map in~\cite{MeyerNest}, the isomorphism
\eqref{intro:fundamentalisomorphism} means that the 
stable Higson corona construction 
represents an \emph{intrinsic} component of the
apparatus surrounding the assembly map. This provides
 a form of converse 
to the Descent Principle.

We conclude by admitting that
we do not know of any geometric obstructions
to~\eqref{intro:canonicalmap} being an isomorphism, save the existence
of a positive scalar curvature metric on a realisation of $BG$ as a
compact spin$^c$\brd{}manifold.  The latter already obstructs
surjectivity of $p_\EG^*$ and \emph{a fortiori} the surjectivity
of~\eqref{intro:canonicalmap} (see \cite{Ros2}). One would therefore
expect there to be an easier obstruction.  Although the existence of
an expanding sequence of graphs embedded in~$\abs{G}$ has negative
consequences for the coarse assembly map for~$\abs{G}$, the arguments
showing this do not appear to apply to the coarse \emph{co}-assembly
map.

\section{Coarse geometry, the stable Higson corona, and the
  coarse co-assembly map}

This article is an expanded and almost completely rewritten version of
the eprint~\cite{EmersonMeyer1}, where we have introduced the stable
Higson corona of a coarse space and the coarse co\brd{}assembly map
and related it to the Dirac dual Dirac method.  In the meantime, we
have discussed the basic properties of these constructions in much
greater detail in~\cite{EmersonMeyer}, within the broader framework of
coarse spaces which need not be groups.  We shall therefore use the
notations and results of~\cite{EmersonMeyer} and concentrate on things
that are special in the group case.  In our analysis, we are going to study
 actions of locally compact groups on
coarse spaces.  This is necessary to prove our descent theorem for
groups with torsion and to define boundary classes.

\subsection{Group actions on coarse spaces}
\label{sec:act_on_coarse}

\begin{defn}  \label{def:act_on_coarse}
  Let~$X$ be a coarse space and let~$G$ be a locally compact group
  that acts continuously on~$X$.  For compatibility with our later
  arguments, we let~$G$ act on the right.  We call the action
  \emph{coarse} if the set of $(xg,yg)$ with $g\in K$, $(x,y)\in E$ is
  an entourage for any compact subset $K\subseteq G$ and any entourage
  $E\subseteq X\times X$.  We say that~$G$ \emph{acts by translations}
  if the set of $(xg,x)$ with $g\in K$, $x\in X$ is an entourage for
  any compact subset $K\subseteq G$.  We say that~$G$ acts
  \emph{isometrically} if any entourage is contained in an
  $G$\nbd{}invariant entourage.
\end{defn}

Actions by translations and isometric actions are coarse.  We usually
require actions to be continuous, isometric and proper.  There are a
few situations where other group actions are also useful.

Let~$G$ be discrete.  Then an action is coarse if and only if the
group acts by coarse maps.  The action is by translations if and only
if these maps are close to the identity.  If the coarse structure
on~$X$ is countably generated, then the action is isometric if and
only if the coarse structure comes from an $G$\nbd{}invariant metric.

Now let $\X=\bigcup X_n$ be a $\sigma$\nbd{}coarse space
(see~\cite{EmersonMeyer} for the definition).  We will only consider
group actions that leave the subsets $X_n$ invariant.  More generally,
it suffices to assume that for all $m\in\N$ there is $n\in\N$ with
$X_m\cdot G\subseteq X_n$; then we can rewrite~$\X$ as $\bigcup
X_n\cdot G$.  We call the action on~$\X$ coarse, isometric, etc., if
the restrictions to~$X_n$ have the corresponding property for all
$n\in\N$.

\subsection{Coarse spaces from groups and proper group
  actions}
\label{sec:coarse_proper}

\begin{thm}  \label{the:proper_act_coarse}
  Let~$G$ be a locally compact group and let~$X$ be a
  \emph{$G$\nbd{}compact}, proper $G$\nbd{}space.  There is a unique
  coarse structure on~$X$ that is compatible with the given topology
  and for which~$G$ acts isometrically.  It is generated by the
  $G$\nbd{}invariant entourages $E_L \defeq \bigcup_{g\in G} Lg\times
  Lg$ for compact $L\subseteq X$.
\end{thm}

We write~$\abs{X}$ for~$X$ equipped with this coarse structure.

\begin{proof}
  It is easy to see that the coarse structure generated by the
  entourages~$E_L$ has the required properties.  Equip~$X$ with any
  coarse structure with the required properties.  Let $L\subseteq X$
  be compact.  Then~$L$ is bounded, so that $L\times L$ is contained
  in a $G$\nbd{}invariant entourage.  Thus~$E_L$ is an entourage.
  Conversely, let $E\subseteq X\times X$ be a $G$\nbd{}invariant
  entourage.  Let $K\subseteq X$ be compact such that $K\cdot G=X$.
  Then $E\cap (K\times X)\subseteq K\times L$ for some bounded and
  hence relatively compact subset~$L$.  We may replace~$L$ by a
  compact subset that contains~$K$.  Since $K\cdot G=X$, the
  $G$\nbd{}invariant entourage~$E$ is determined by $E\cap (K\times
  X)$.  We obtain $E\subseteq E_L$.  Hence the coarse structure is
  equal to the one defined by the entourages~$E_L$.
\end{proof}

In particular, we let~$\abs{G}$ be the group~$G$ itself equipped with
the action of~$G$ by right multiplication and with the unique coarse
structure for which this action is isometric.

Let~$X$ be a $G$\nbd{}compact proper $G$\nbd{}space.  For any $x\in
X$, the map $\abs{G}\to\abs{X}$, $g\mapsto g\cdot x$, is a coarse
equivalence.  These maps for different points in~$X$ are close.

Now let~$X$ be a proper $G$\nbd{}space that is not necessarily
$G$\nbd{}compact.  We only require~$X$ to be a union of an increasing
sequence $(X_n)_{n\in\N}$ of $G$\nbd{}compact subspaces.  We
implicitly require the~$X_n$ to be $G$\nbd{}invariant and closed.
Even if~$X$ is not locally compact, the spaces~$X_n$ are necessarily
locally compact in the subspace topology.  Thus we can turn them into
coarse spaces by the above prescription.  The maps
$\abs{X_m}\to\abs{X_n}$ for $m\le n$ are coarse equivalences because
orbit maps $\abs{G}\to\abs{X_m}$ are coarse equivalences.  Hence we
have turned~$X$ into a $\sigma$\nbd{}coarse space in the sense
of~\cite{EmersonMeyer}.  We write~$\abs{X}$ for this
$\sigma$\nbd{}coarse space.

Of course, the above construction is natural, that is, a
$G$\nbd{}equivariant continuous map $f\colon X\to Y$ induces a coarse
continuous map $\abs{X}\to\abs{Y}$.

\subsection{The coarse category of coarse spaces}
\label{sec:coarse}

From now on, we require locally compact groups and topological spaces
to be second countable and coarse structures to be countably
generated.  We let~$\coa$ be the category of coarse spaces, whose
objects are the coarse spaces with second countable topology and
countably generated coarse structure and whose morphisms are the
continuous coarse maps.  Let~$H$ be a second countable locally compact
group.  We let~$\coa_H$ be the category of coarse spaces as above,
equipped with a continuous, proper, and isometric action of~$H$.  The
morphisms in~$\coa_H$ are the $H$\nbd{}equivariant coarse continuous
maps.

Let $X\in\coa_H$ and let $Y\subseteq X$ be an $H$\nbd{}invariant
closed subset.  Then we give~$Y$ the subspace coarse structure, so
that $Y\in\coa_H$.  We say that~$Y$ is \emph{coarsely dense} if there
is an entourage $E\subseteq X\times X$ such that for any $x\in X$
there is $(x,y)\in E$ with $y\in Y$.  If~$Y$ is discrete, we call~$Y$
a \emph{discretisation} of~$X$.  A discretisation exists if and only
if the orbits of the $H$\nbd{}action are discrete.

If~$H$ is finite, we want all finite $H$\nbd{}spaces to be coarsely
equivalent.  However, there is no $H$\nbd{}equivariant map from the
one point space with trivial action to~$H$.  To overcome this problem,
we relax $H$\nbd{}equivariance as follows:

\begin{defn}  \label{def:approx_equivariant}
  Let $X,Y\in\coa_H$ and let $f\colon X\to Y$ be a coarse Borel
  map.  Define $f^h(x)\defeq f(xh)h^{-1}$ for $x\in X$, $h\in H$.  We
  call~$f$ \emph{almost $H$\nbd{}equivariant} if
  $$
  \bigl\{\bigl(f^h(x),f(x)\bigr) \bigm| x\in X,\ h\in H\bigr\}
  $$
  is an entourage.
\end{defn}

\begin{defn}  \label{def:coarse_category}
  We let $\cocoa_H$ be the category with the same objects as $\coa_H$
  and whose morphisms are equivalence classes of almost
  $H$\nbd{}equivariant coarse Borel maps, where two maps are
  identified if they are close.  We call $\cocoa_H$ the \emph{coarse
    category of coarse $H$\nbd{}spaces}.  A morphism $f\colon X\to Y$
  is called a \emph{coarse equivalence} if it is an isomorphism in
  $\cocoa_H$.
\end{defn}

The following lemma is easy to check:

\begin{lemma}
  The embedding of a coarsely dense $H$\nbd{}invariant subspace is a
  coarse equivalence.
\end{lemma}

\subsection{An equivariant version of the Rips complex}
\label{sec:equivariant_Rips}

The Rips complex can be constructed most easily for discrete coarse
spaces.  In~\cite{EmersonMeyer}, we have defined it for non-discrete
spaces by choosing a discretisation.  When we work equivariantly, this
does not always exist.  To get a more natural construction of
$\Rips(X)$ for non-discrete~$X$, we adapt the locally compact model
for $\EG$ for a locally compact group~$G$ due to Gennadi Kasparov and
Georges Skandalis (\cite{KasparovSkandalis2}).

Let $X\in\coa_H$.  Let $\Rips(X)$ be the set of all positive Borel
measures~$\mu$ on~$X$ with total volume $1/2<\norm{\mu}_1\le1$.  We
equip (subsets of) $\Rips(X)$ with the weak topology from the pairing
with $C_0(X)$.  In this topology, $\Rips(X)$ is locally compact and
second countable.  Let $E\subseteq X\times X$ be an entourage.  A
measurable subset $S\subseteq X$ is called \emph{$E$\nbd{}bounded} if
$S\times S \subseteq E$.  Given a closed entourage~$E$ and $t>1/2$, we
let
$$
\Subrips_{E,t}(X)\defeq \{\mu\in\Rips(X) \mid
  \text{$\mu(S)\ge t$ for some $E$-bounded set $S$}
\}.
$$
Since~$X$ is countably generated, there is an increasing sequence
of closed entourages $(E_n)$ that dominates any other entourage.
Write $\Subrips_n(X)\defeq\Subrips_{E_n,1/2+1/n}$.  Then $\Rips(X) =
\bigcup \Subrips_n(X)$.

We claim that $\Subrips_{E,t}(X)$ is a weakly closed subset of
$\Rips(X)$.  Let $(\mu_n)$ be a weakly convergent sequence in
$\Subrips_{E,t}(X)$ that converges towards some $\mu\in\Rips(X)$.
Choose $E$\nbd{}bounded subsets $S_n$ for $n\in\N$ such that
$\mu_n(S_n)\ge t$ and choose a compact subset~$S$ with $\mu(S)\ge t$.
Then $\mu_n(S)>1/2$ and hence $S_n\cap S\neq\emptyset$ for almost
all~$n$.  Since the subsets~$S_n$ are compact and $\bigcap
S_n\neq\emptyset$, the subset $\bigcup S_n$ is relatively compact.
The set of compact subsets of its closure is compact in the Hausdorff
metric.  Hence we can find a subsequence of $(S_n)$ that converges
towards some compact subset $S_\infty$.  We may assume that $(S_n)$
itself converges.  The convergence $\lim \mu_n=\mu$ implies
$\mu(S_\infty)\ge t$.  The convergence $S_n\to S$ implies
$S_\infty\times S_\infty\subseteq E$ because~$E$ is closed.  Thus
$\mu\in\Subrips_{E,t}$, so that $\Subrips_{E,t}(X)$ is closed as
asserted.  It follows that $\Subrips_{E,t}(X)$ is locally compact in
the subspace topology.

We also let
\begin{multline*}
\Subrips^2_{E,t}(X)\defeq
\{(\mu,\nu)\in\Rips(X)\times\Rips(X) \mid
\\ \text{$\mu(S)\ge t$ and $\nu(S)\ge t$ for some $E$-bounded subset
    $S$}
\}.
\end{multline*}
These subsets of $\Rips(X)\times\Rips(X)$ define a coarse structure on
$\Rips(X)$.  One checks easily that the restriction of this coarse structure 
to
$\Subrips_{E',t'}(X)$ is compatible with the topology for all closed
entourages~$E'$ and all $t'>1/2$.  Moreover,
$\Subrips_{E,t}(X)\subseteq\Subrips_{E',t'}(X)$ is coarsely dense if
$E\subseteq E'$ and $t\ge t'$.  Thus $\Rips(X)=\bigcup \Subrips_n(X)$
is a $\sigma$\nbd{}coarse space.

The action of~$H$ on~$X$ induces an action on $\Rips(X)$, pushing
forward measures.  If~$E$ is an $H$\nbd{}invariant
entourage, then $\Subrips_{E,t}(X)$ is $H$\nbd{}invariant and the
restriction of the action to $\Subrips_{E,t}(X)$ is still isometric,
continuous and proper.  That is, $\Subrips_{E,t}(X)\in\coa_H$ for all
$H$\nbd{}invariant entourages~$E$ and all $t>1/2$.

There is a canonical map $j_X\colon X\to\Rips(X)$, sending $x\in X$ to
the Dirac measure at~$X$.  The construction of $\Rips(X)$ is natural:
a continuous coarse map $f\colon X\to X'$ induces a continuous coarse
map $\Rips(f)\colon \Rips(X)\to\Rips(X')$ which pushes forward
measures along~$f$.  We have $\Rips(f)\circ j_X=j_{X'}\circ f$, that
is, $j_X\colon X\to\Rips(X)$ is a natural transformation.

The above construction has another property that is useful for
technical purposes.  For any pair $(E,t)$, there is $n\in\N$ such that
$\Subrips_n(X)$ is a neighbourhood of $\Subrips_{E,t}(X)$ in
$\Rips(X)$.  Passing to a subsequence, we can achieve that
$\Subrips_{n+1}(X)$ is a neighbourhood of $\Subrips_n(X)$ for all
$n\in\N$.  Hence there is a partition of unity $(\phi_n)$ on
$\Rips(X)$ with $\supp\phi_n\subseteq \Subrips_{n+1}(X)\setminus
\Subrips_{n-1}(X)$.  Even more, we can choose the functions $(\phi_n)$
to be $H$\nbd{}equivariant.  These partitions of unity are useful to
formulate the universal property of $\Rips(X)$.  We call a
$\sigma$\nbd{}coarse $H$\nbd{}space $\X=\bigcup X_n$
\emph{partitionable} if it carries such a partition of unity.  This
holds if and only if $X_{n+1}$ is a neighbourhood of~$X_n$
in $X_{n+2}$ for all $n\in\N$.

We let $\sigma\coa_H$ be the category of partitionable
$\sigma$\nbd{}coarse $H$\nbd{}spaces.

Let $X\in\coa_H$.  We give $X\times [0,1]$ the product topology and
let $E\subseteq (X\times [0,1])^2$ be an entourage if its image in
$X\times X$ is one.  Then $X\times[0,1]\in\coa_H$.  The evaluation
maps $X\cong X\times\{t\}\subseteq X\times [0,1]$ and the projection
$X\times [0,1]\to X$ are coarse equivalences.  A \emph{coarse
  homotopy} between two morphisms $f,g\colon X\to Y$ is an
$H$\nbd{}equivariant coarse continuous map $X\times[0,1]\to Y$.  This
generates the equivalence relation of \emph{coarse homotopy} on the
space of morphisms $\coa_H(X,Y)$, which in turn generates a notion of
\emph{coarse homotopy equivalence}.  If two maps are coarsely
homotopic, then they are both close and homotopic (as maps of
topological spaces).  Thus a coarse homotopy equivalence is
simultaneously a coarse equivalence and a homotopy equivalence.  It is
evident how to extend these notions to $\sigma$\nbd{}coarse spaces.

\begin{lemma}  \label{lem:Rips_universal}
  Let $X\in\coa_H$.  Then $\Rips(X)\in\sigma\coa_H$, and $j_X\colon
  X\to\Rips(X)$ is a coarse equivalence.  Let $\Y\in\sigma\coa_H$ and
  let $f\colon X\to\Y$ be a coarse equivalence.  Then there is
  $h\in\sigma\coa_H\bigl(\Y,\Rips(X)\bigr)$ such that $h\circ f$ is
  coarsely homotopy equivalent to~$j_X$, and the map~$h$ is unique up
  to coarse homotopy equivalence.  This universal property determines
  $\Rips(X)$ uniquely up to coarse homotopy equivalence.
\end{lemma}

\begin{proof}
  It is straightforward to see that $\Rips(X)$ has the required
  properties.  The map~$j_X$ is a coarse equivalence because it is an
  embedding with coarsely dense range.  If $f_0,f_1\colon
  \Y\to\Rips(X)$ are two close maps, then we can join them by the
  affine homotopy $(1-t)f_0+tf_1$.  One checks easily that this
  homotopy is coarse.  Thus close coarse maps into $\Rips(X)$ are
  coarsely homotopic.  Therefore, the proof will be finished if we
  construct a coarse continuous map $h\colon \Y\to\Rips(X)$ such that
  $h\circ f$ is close to the identity map.

  Write $\Y=\bigcup Y_n$ and suppose that we have found maps
  $h_n\colon Y_n\to\Rips(X)$ with the required properties.  Since~$\Y$
  is partitionable, we obtain a certain partition of unity $(\phi_n)$.
  Since $\Rips(X)$ is convex, we can define $h\defeq \sum \phi_n
  h_{n+1}$.  This map has the required properties.  Thus it remains to
  construct maps $h_n\colon Y_n\to\Rips(X)$.  Since~$f$ is a coarse
  equivalence, there is an almost $H$\nbd{}equivariant Borel map
  $g\colon Y_n\to X$ such that $gf$ is close to the identity.  Choose
  a uniformly bounded open covering of~$Y_n$ and some subordinate
  partition of unity $(\psi_k)_{k\in\N}$, and choose $x_k\in X$ close
  to $g(\supp\psi_k)$.  Define $g'(y)\defeq \sum_{k\in\N} \psi_k(y)
  \delta_{x_k}$, where $\delta_{x_k}$ denotes the Dirac measure
  at~$x_k$.  This is a continuous map $g'\colon Y_n\to\Rips(X)$ that
  is close to~$g$.  Hence~$g'$ is almost $H$\nbd{}equivariant and
  coarse and $g'\circ f$ is close to~$j_X$.  To achieve exact
  $H$\nbd{}equivariance, we choose a cut off function~$\Psi$ on~$Y_n$.
  This exists because the action is proper.  We obtain a map $h_n\colon
  Y_n\to\Rips(X)$ with the required properties if we define
  $$
  \langle h_n(y), \alpha\rangle
  \defeq \int_H \langle g'(ys)\cdot s^{-1},\alpha \rangle \,
  \Psi(ys) \,ds
  $$
  for all $y\in Y_n$, $\alpha\in C_0(X)$.
\end{proof}

\begin{lemma}  \label{lem:Rips_EG}
  Let~$G$ be a second countable locally compact group and let $\EG$ be
  any universal proper $G$\nbd{}space.  Then $\abs\EG$ is
  $G$\nbd{}equivariantly coarsely homotopy equivalent to
  $\Rips(\abs{G})$.
\end{lemma}

\begin{proof}
  If we choose the special model for $\EG$ constructed
  in~\cite{KasparovSkandalis2}, then we have $\Rips(\abs{G})=\abs\EG$.
  Since $\EG$ is determined uniquely up to $G$\nbd{}equivariant
  homotopy equivalence and since the construction $X\mapsto\abs{X}$ is
  natural, any two models for $\abs\EG$ are $G$\nbd{}equivariantly
  coarsely homotopy equivalent.
\end{proof}

\begin{lemma}  \label{lem:coarse_maps_Rips}
  The set $\cocoa_H(X,Y)$ of morphisms $X\to Y$ in the coarse category
  is naturally isomorphic to the set of coarse homotopy classes of
  $H$\nbd{}equivariant coarse continuous maps $\Rips(X)\to\Rips(Y)$.
\end{lemma}

\begin{proof}
  As in the proof of Lemma~\ref{lem:Rips_universal} one shows that any
  almost equivariant coarse Borel map $X\to Y$ is close to an
  $H$\nbd{}equivariant coarse continuous map $X\to\Rips(Y)$.
  Moreover, close maps into $\Rips(Y)$ are coarsely homotopic.  Thus
  $\cocoa_H(X,Y)$ is in bijection with coarse homotopy classes of
  coarse continuous maps $X\to\Rips(Y)$.  Write $\Rips(Y)=\bigcup
  \Subrips_n(Y)$.  Any map $X\to\Rips(Y)$ is a map into
  $\Subrips_n(Y)$ for some~$n$ and hence induces a map
  $\Rips(X)\to\Rips(\Subrips_n(Y))$.  Since $Y\to\Subrips_n(Y)$ is a
  coarse equivalence, the induced map
  $\Rips(Y)\to\Rips(\Subrips_n(Y))$ is a coarse homotopy equivalence
  and hence has an inverse up to coarse homotopy.  It follows that any
  morphism $X\to\Rips(Y)$ is close to the restriction of a morphism
  $\Rips(X)\to\Rips(Y)$.  As above, this implies that $\cocoa_H(X,Y)$
  is in bijection with coarse homotopy classes of morphisms
  $\Rips(X)\to\Rips(Y)$.
\end{proof}

Thus the map $X\to\Rips(X)$ plays the role of an injective resolution.
If we are given any functor~$F$ on $\coa_H$, we obtain a functor that
descends to $\cocoa_H$ by applying~$F$ to $\Rips(X)$ instead of~$X$.

For discrete spaces, the space $\Rips(X)$ is usually constructed using
probability measures with compact support.  We can also do this in
general.  An elementary argument shows that the two versions for
$\Rips(X)$ are coarsely homotopy equivalent, so that it makes no
difference which one we use.  However, the model for $\Rips(X)$ with
compactly supported probability measures is not partitionable.  This
makes it more difficult to formulate the universal property.  In any
case, for explicit computations, one tends to look for smaller models
for $\Rips(X)$.  For instance, as is well-known (and proved
in~\cite{EmersonMeyer}), we can use~$X$ itself if~$X$ has bounded
geometry and is uniformly contractible.

\subsection{Coarse K-theory}
\label{sec:coarse_Ktheory}

Let $\Y=\bigcup Y_n$ be any $\sigma$\nbd{}coarse space.  For our
purposes, the appropriate algebra of functions on~$Y$ is the
$\sigma$\nbd{}$C^*$\brd{}algebra $C_0(\Y)\defeq \varprojlim C_0(Y_n)$,
which consists of all functions $f\colon \Y\to\C$ for which
$f|_{Y_n}\in C_0(Y_n)$ for all $n\in\N$ (see
Lemma~\ref{lem:RKKG_as_Ktheory}).  An action of~$H$ on~$\Y$ turns
$C_0(\Y)$ into a projective system of $H$\nbd{}$C^*$\brd{}algebras,
which we call a \emph{$\sigma$\nbd{}$H$\brd{}$C^*$\brd{}algebra}.  We
define the crossed product for a
$\sigma$\nbd{}$H$\brd{}$C^*$\brd{}algebra $\varprojlim A_m$ in the
evident way as $(\varprojlim A_m)\cross H\defeq \varprojlim A_m\cross
H$.  This yields a $\sigma$\nbd{}$C^*$\brd{}algebra.

We define the \emph{coarse $\K$\nbd{}theory} $\KX^*(X)$ of a coarse
space~$X$ as the $\K$\nbd{}theory of the
$\sigma$\nbd{}$C^*$\brd{}algebra $C_0\bigl(\Rips(X)\bigr)$.  We define
the \emph{$H$\nbd{}equivariant coarse $\K$\nbd{}theory} $\KX^*_H(X)$
of $X\in\coa_H$ as the $\K$\nbd{}theory of the
$\sigma$\nbd{}$C^*$\brd{}algebra $C_0(\Rips(X)\cross H)$.  We may also
tensor $C_0(\Rips(X))$ with an $H$\nbd{}$C^*$\brd{}algebra~$D$.  This
yields the $H$\nbd{}equivariant coarse $\K$\nbd{}theory $\KX^*_H(X,D)$
of~$X$ with coefficients~$D$.  All the above groups are evidently
functorial for coarse continuous maps $\Rips(X)\to\Rips(Y)$, and
(coarsely) homotopy equivalent maps induce the same map.  Hence we
obtain functors on the coarse category $\cocoa_H$ by
Lemma~\ref{lem:coarse_maps_Rips}.

Let~$G$ be a locally compact group and let $H\subseteq G$ be a closed
subgroup.  Then
$$
\KX^*_H(\abs{G},D) \cong
\K_*\bigl((C_0(\abs\EG)\otimes D)\cross H\bigr)
$$
by Lemma~\ref{lem:Rips_EG}.  \emph{We warn the reader that
  $C_0(\abs{\EG},D)\neq C_0(\EG,D)$ unless $\EG$ is $G$\nbd{}compact}:
elements of $C_0(\abs{\EG},D)$ are possibly unbounded continuous
functions $\EG\to D$, only their restrictions to $G$\nbd{}compact
subsets vanish at~$\infty$.

\subsection{The stable Higson corona}
\label{sec:bvr}

Let~$X$ be a coarse space and~$D$ a $C^*$\nbd{}algebra.  Let
$\Mult(D\otimes\Comp)$ be the multiplier algebra of $D\otimes\Comp$.
We identify it with the $C^*$\nbd{}algebra of adjointable operators on
the Hilbert $D$\nbd{}module $D\otimes\ell^2(\N)$.  Let
$\bbbarred(X,D)$ be the $C^*$\nbd{}algebra of norm continuous, bounded
functions $f\colon X\to\Mult(D\otimes\Comp)$ for which $f(x)-f(y)\in
D\otimes\Comp$ for all $x,y\in X$.  We say that a function
$f\in\bbbarred(X,D)$ has \emph{vanishing variation} if the function
$E\ni (x,y)\mapsto \norm{f(x)-f(y)}$ vanishes at~$\infty$ for any
closed entourage $E\subseteq X\times X$.  Let $\vvr(X,D) \subseteq
\bbbarred(X,D)$ be the subalgebra of vanishing variation functions.
Let
\begin{align*}
  \bvr(X,D) &\defeq \vvr(X,D)/C_0(X,D\otimes\Comp),
  \\
  \bbred(X,D) &\defeq \bbbarred(X,D)/C_0(X,D\otimes\Comp).
\end{align*}
It is clear that $X\mapsto \vvr(X,D)$ and $X\mapsto \bvr(X,D)$ are
contravariant functors in~$X$ for continuous coarse maps.  If $D=\C$,
we drop it from our notation and write $\vvr(X)$ and $\bvr(X)$.  We
refer to $\vvr(X)$ and $\bvr(X)$ as the \emph{reduced stable Higson
  compactification} and the \emph{reduced stable Higson corona}
of~$X$, respectively.

\begin{prop}[\cite{EmersonMeyer}]  \label{pro:bvr_coarse}
  The functor $X\mapsto\bvr(X,D)$ descends to a functor on the coarse
  category of coarse spaces.  That is, close maps $f,f'\colon X\to X'$
  induce the same map $\bvr(X',D)\to\bvr(X,D)$ and a coarse
  equivalence $X\to X'$ induces an isomorphism
  $\bvr(X',D)\cong\bvr(X,D)$.
\end{prop}

Now let $\X=\bigcup X_n$ be a $\sigma$\nbd{}coarse space.  We have
already defined $C_0(\X,D)\defeq \varprojlim C_0(X_n,D)$ above.  Let
\begin{displaymath}
  \vvr(\X,D) \defeq \varprojlim \vvr(X_n,D),
  \qquad
  \bvr(\X,D) \defeq \varprojlim \bvr(X_n,D).
\end{displaymath}
Equivalently, $\vvr(\X,D)$ is the $\sigma$\nbd{}$C^*$\brd{}algebra of
all functions $f\colon \X\to\Mult(D\otimes\Comp)$ for which
$f|_{X_n}\in\vvr(X_n,D)$ for all $n\in\N$.  Since the maps $X_m\to
X_n$ for $m\le n$ are coarse equivalences,
Proposition~\ref{pro:bvr_coarse} yields that
$\bvr(\X,D)\cong\bvr(X_n,D)$ for all $n\in\N$, so that $\bvr(\X,D)$ is
still a $C^*$\nbd{}algebra.  It is observed in~\cite{EmersonMeyer}
that the obvious maps give rise to an extension of
$\sigma$\nbd{}$C^*$\brd{}algebras
$$
0\to C_0(\X,D\otimes\Comp)\to\vvr(\X,D)\to\bvr(\X,D)\to0.
$$

Now let a locally compact group~$H$ act continuously and coarsely on a
coarse space~$X$ and let~$D$ be an $H$\nbd{}$C^*$\brd{}algebra.  Let
$\Comp_H\defeq \Comp(\ell^2\N\otimes L^2 H)$ and let~$H$ act on
$D\otimes\Comp_H$ and $\Mult(D\otimes\Comp_H)$ in the obvious way.  We
let $\bbbarred_H(X,D)$ be the $H$\nbd{}continuous subspace of the
$C^*$\nbd{}algebra of norm continuous, bounded functions $f\colon
X\to\Mult(D\otimes\Comp_H)$ for which $f(x)-f(y)\in D\otimes\Comp_H$
for all $x,y\in X$.  The group~$H$ acts on $\bbbarred_H(X,D)$ by
$(h\cdot f)(x)\defeq h\cdot \bigl(f(xh)\bigr)$.  As above, we let
$\vvr_H(X,D)$ be the subspace of vanishing variation functions in
$\bbbarred_H(X,D)$.  This subalgebra is $H$\nbd{}invariant.  It
contains $C_0(X,D\otimes\Comp_H)$ as an $H$\nbd{}invariant ideal.
Thus $\bvr_H(X,D)\defeq \vvr_H(X,D)/C_0(X,D\otimes\Comp_H)$ is an
$H$\nbd{}$C^*$\brd{}algebra.  The same holds for $\bbred_H(X,D)\defeq
\bbbarred_H(X,D)/C_0(X,D\otimes\Comp_H)$.

It is necessary to restrict attention to $H$\nbd{}continuous elements
in order to form the crossed product $\bvr_H(X,D)\cross H$.  We warn
the reader that, although a vanishing variation function $\abs{G}\to
D$ is automatically $G$\nbd{}continuous for the action by left
multiplication, which is by translations, this implies nothing about
the action by right multiplication, which is the one we use.

The following is the central definition of this paper.

\begin{defn}  \label{def:equivariant_corona_Ktheory}
  We refer to $\Ktop_*\bigl(H,\bvr_H(X,D)\bigr)$ as the
  \emph{$H$\nbd{}equivariant boundary $\K$\nbd{}theory of~$X$ with
    coefficients~$D$}.
\end{defn}

We are particularly interested in the case where~$G$ is a discrete
group, $D=\C$ and $X=\abs{G}$, and where~$H$ is a subgroup of~$G$
acting by right multiplication on~$G$.  When~$H$ is the trivial group,
we obtain $\K_*\bigl(\bvr(\abs{G})\bigr)$, the $\K$\nbd{}theory of the
stable Higson corona of~$G$, which is used for descent.  The
equivariant version $\Ktop_*\bigl(G,\bvr(\abs{G})\bigr)$ gives the
boundary classes in $\RKK^G_*(\EG;\C,\C)$.

\subsection{The coarse co-assembly map}
\label{sec:coarse_coassembly}

Let $X\in\coa_H$ and let~$D$ be an $H$\nbd{}$C^*$\brd{}algebra.  Form
the Rips complex $\Rips(X)$ as above.  Let $\Dirac\in\KK^G(\ADir,\C)$
be a Dirac morphism for~$G$ as in~\cite{MeyerNest}.  Recall that we
identify $\Ktop_*(H,D)\defeq \K_*\bigl((D\otimes\ADir)\cross H\bigr)$.
We also use this definition for
$\sigma$\nbd{}$H$\brd{}$C^*$\brd{}algebras.  Since the action of~$H$
on $\Rips(X)$ is proper, $C_0(\Rips(X),D)$ is an inverse system of
proper $G$\nbd{}$C^*$\brd{}algebras $(B_m)$.  These belong to
$\gen{\CI}$ by \cite{MeyerNest}*{Corollary 7.3}.  Hence the maps
$\Dirac_*\colon (B_m\otimes\ADir)\cross H\to B_m\cross H$ are
isomorphisms on $\K$\nbd{}theory.  It follows from the Milnor
$\varprojlim^1$\brd{}sequence of~\cite{Phillips} that the map
$\K_*(\varprojlim B_m\otimes\ADir)\cross H)\to \K_*(\varprojlim
B_m\cross H)$ is an isomorphism as well.  Thus
\begin{displaymath}
  \K_*\bigl((C_0(\Rips(X),D)\otimes\ADir)\cross H\bigr)
  \cong \KX^*_H(X,\ADir\otimes D)
  \cong \KX^*_H(X,D).
\end{displaymath}

The coarse equivalence $j_X\colon X\to\Rips(X)$ induces an isomorphism
$$
\bvr(\Rips(X),D)\cong\bvr(X,D).
$$
We have a canonical extension of
$\sigma$\nbd{}$H$\brd{}$C^*$\brd{}algebras
\begin{displaymath}
  0\to C_0(\Rips(X),D\otimes\Comp_H)
  \to \vvr_H(\Rips(X),D)
  \to \bvr_H(\Rips(X),D)
  \to 0.
\end{displaymath}
Taking (maximal) tensor products with the source of the Dirac
morphism~$\ADir$ and (full) crossed products with~$H$, we obtain an
extension of $\sigma$\nbd{}$C^*$\brd{}algebras
\begin{multline}  \label{eq:extension_vvr_bvr_crossed}
  0\to (C_0(\Rips(X),D\otimes\Comp_H)\otimes\ADir) \cross H
  \to (\vvr_H(\Rips(X),D)\otimes\ADir)\cross H
  \\ \to (\bvr_H(\Rips(X),D)\otimes\ADir)\cross H
  \to 0.
\end{multline}
We have canonical isomorphisms
\begin{equation}  \label{somecanonicalisomorphisms}
\begin{aligned}
  \K_*\bigl((\bvr_H(\Rips(X),D)\otimes\ADir)\cross H\bigr)
  &\cong \Ktop_*\bigl(H,\bvr_H(X,D)\bigr),
  \\
  \K_*\bigl((C_0(\Rips(X),D\otimes\Comp_H)\otimes\ADir) \cross H\bigr)
  &\cong \KX^*_H(X,D).
\end{aligned}
\end{equation}
Via these isomorphisms, the $\K$\nbd{}theory boundary map
for~\eqref{eq:extension_vvr_bvr_crossed} is equivalent to a map
\begin{equation}  \label{defnofcoarsecoassemblymap}
\mu_{X,H,D}^*\colon
\Ktop_{*+1}\bigl(H,\bvr_H(X,D)\bigr) \to
\KX^*_H(X,D).
\end{equation}

\begin{defn}  \label{def:coarsecoassemblymap}
  We call~\eqref{defnofcoarsecoassemblymap} the
  \emph{$H$\nbd{}equivariant coarse co\brd{}assembly map for~$X$ with
    coefficients~$D$}.
\end{defn}

We obtain an equivalent map if we replace $\Rips(X)$ by a coarsely
homotopy equivalent object of $\sigma\coa_H$.  The map $\mu_{X,H,D}^*$
is natural for almost equivariant coarse Borel maps by
Lemma~\ref{lem:coarse_maps_Rips}; that is, any such map gives rise to
a commuting square diagram.  The map $\mu_{X,H,D}^*$ is an isomorphism
if and only if $\Ktop_*\bigl(H,\vvr_H(\Rips(X),D)\bigr) = 0$.

If~$G$ is a locally compact group and $\EG$ is a universal proper
$G$\nbd{}space, then $\abs\EG$ is coarsely homotopy equivalent to
$\Rips(\abs{G})$ by Lemma~\ref{lem:Rips_EG}.  Hence we may use
$\abs\EG$ instead of $\Rips(\abs{G})$ to construct
$\mu_{\abs{G},H,D}^*$ for any closed subgroup $H\subseteq G$, acting
on~$\abs{G}$ by right multiplication.

For instance, let~$G$ be a discrete group.  Set $D=\C$ and
$X=\abs{G}$.  We may then allow~$H$ to run through the subgroups
of~$G$.  If~$H$ is the trivial subgroup we obtain the ordinary coarse
co-assembly map for~$\abs{G}$:
$$
\mu_{\abs{G}}^*\colon
\K_{*+1}\bigl( \bvr(\abs{G})\bigr) \to
\KX^*(\abs{G}) \cong \K_*\bigl( C_0 (\abs{\EG})\bigr).
$$
If $H=G$, we obtain a map
$$
\mu_{\abs{G},G}^*\colon
\Ktop_{*+1}\bigl(G,\bvr(\abs{G})\bigr) \to
\KX^*_G(\abs{G}) \cong \K_*\bigl( C_0 (\abs{\EG})\cross G\bigr).
$$

\subsection{The co-assembly map as a forgetful map}
\label{sec:mu_forget}

Let~$X$ be a uniformly contractible coarse space with bounded geometry
(without group action).  It is shown in~\cite{EmersonMeyer} that the
natural map $X\to\Rips(X)$ is a coarse homotopy equivalence in this
case.  Hence $\KX^*(X,D)\cong\K_*(C_0(X,D))$, and the coarse
co\brd{}assembly map is equivalent to the connecting map for the
extension
$$
0 \to C_0(X,D\otimes\Comp) \to \vvr(X,D) \to \bvr(X,D) \to 0.
$$
We can express the latter as a \emph{forget-control map} as
follows.  By definition, we have $ \vvr(X,D)\subseteq\bbbarred(X,D)$.
This yields a morphism of $C^*$\nbd{}extensions
$$
\xymatrix{
  0 \ar[r] &
  C_0(X,D\otimes\Comp) \ar[r] \ar@{=}[d] &
  \vvr(X,D) \ar[r] \ar[d] &
  \bvr(X,D) \ar[r] \ar[d]^{j} &
  0 \\
  0 \ar[r] &
  C_0(X,D\otimes\Comp) \ar[r] &
  \bbbarred(X,D) \ar[r] &
  \bbred(X,D) \ar[r] &
  0.
}
$$

\begin{prop}  \label{pro:mu_forgetful}
  If~$X$ is a uniformly contractible coarse space of bounded geometry,
  then the $\K$\nbd{}theory connecting map $\partial'\colon
  \K_{*+1}\bigl(\bbred(X,D)\bigr) \to \K_*\bigl(C_0(X,D)\bigr)$ is an
  isomorphism.  Hence $\mu_{X,D}^*$ is equivalent to the map
  $$
  j_*\colon \K_{*+1}\bigl(\bvr(X,D)\bigr) \to
  \K_{*+1}\bigl(\bbred(X,D)\bigr).
  $$
\end{prop}

To prove Proposition~\ref{pro:mu_forgetful}, it suffices to show that
the $C^*$\nbd{}algebra $\bbbarred(X,D)$ has vanishing $\K$\nbd{}theory
whenever~$X$ is uniformly contractible.  Actually, this will be the
case under the weaker assumption that~$X$ is contractible.

Let $\bbbar(X,D)$ denote the $C^*$-algebra of bounded, norm continuous
maps $X\to D\otimes\Comp$.  There is a canonical embedding
$D\otimes\Comp \to\bbbar(X,D)$.  We have:

\begin{lemma}  \label{identificationofreducedktheoryofbbalgebras}
  The inclusion $\bbbar(X,D)\to\bbbarred(X,D)$ induces a canonical
  isomorphism
  $$
  \K_*\bigl(\bbbar(X,D)\bigr)/\K_*(D)\cong \K_*\bigl(\bbbarred(X,D)\bigr).
  $$
\end{lemma}

\begin{proof}
  Let~$p$ be the composition $\bbbarred(X,D) \to\Bound(D\otimes\Comp)
  \to \Calk(D\otimes\Comp)$, where the first map is evaluation at some
  point $\star\in X$, the second map is the quotient map, and where
  $\Calk(D\otimes \Comp)$ denotes
  $\Bound(D\otimes\Comp)/D\otimes\Comp$.  The induced map on
  $\K$\nbd{}theory vanishes because it factorises through the
  $\K$\nbd{}theory of $\Bound(D\otimes\Comp) $, which is zero.  The
  kernel of~$p$ is evidently $\bbbar(X,D)$.  By considering the
  associated six-term exact sequence and making the identification
  $\K_{*+1}(\Calk(D\otimes \Comp))\cong\K_*(D)$, one obtains short
  exact sequences
  $$
  0\to \K_*(D) \to\K_*\bigl(\bbbar(X,D)\bigr)
  \to\K_*\bigl(\bbbarred(X,D)\bigr) \to 0.
  $$
  It is easily checked that the first map is the canonical
  inclusion.
\end{proof}

\begin{lemma}  \label{homotopyinvarianceofbbbar}
  For any~$D$ there is a canonical isomorphism
  $$
  \K_{*+1}\bigl(\bbbar(X,D)\bigr) \cong \RKK_*(X;\C,D).
  $$
\end{lemma}

\begin{proof}
  We claim that $\Mult\bigl(C_0(X,D\otimes\Comp)\bigr) \cong
  \Mult\bigl(\bbbar(X,D)\bigr)$.  The multipliers of
  $C_0(X,D\otimes\Comp)$ are the bounded, strictly continuous
  functions $X\to\Mult(D\otimes\Comp)$.  These also act as multipliers
  on $\bbbar(X,D)$.  The converse is clear because
  $C_0(X,D\otimes\Comp)$ is an ideal in $\bbbar(X,D)$.  It follows
  that the $\K$\nbd{}theory of $\Mult\bigl(\bbbar(X,D)\bigr)$
  vanishes.  Thus the $\K$-theory of $\bbbar(X,D)$ is isomorphic to
  the $\K$-theory of $\Mult\bigl( \bbbar(X,D)\bigr)/\bbbar(X,D)$.  It
  remains to identify the $\K$\nbd{}theory of the latter with
  $\RKK_*(X;\C,D)$.  We sketch the proof for $*=0$.  Cycles for
  $\RKK_0(X;\C,D)$ are $F\in \Mult\bigl(C_0(X,D\otimes\Comp)\bigr)$
  for which $M_\phi\cdot (FF^*-1)$ and $M_\phi\cdot (F^*F-1)$ lie in
  $C_0(X,D\otimes\Comp)$ for all $\phi\in C_0(X)$.  Equivalently,
  $FF^*-1$ and $F^*F-1$ belong to $\bbbar(X,D)$.

  Two cycles differ by a compact perturbation if and only if their
  difference belongs to $\bbbar(X,D)$.  Thus equivalence classes of
  cycles for $\RKK_0(X;\C,D)$ up to compact perturbation are the same
  as unitaries in $\Mult\bigl(\bbbar(X,D)\bigr)/\bbbar(X,D)$.
  Moreover, operator homotopy of cycles corresponds to homotopy of
  unitaries, and degenerate cycles correspond to unitaries in
  $\Mult\bigl(\bbbarred(X,D)\bigr)$.  These observations together with
  the homotopy invariance of Kasparov theory yields the assertion.  A
  similar argument also occurs in the proof of
  Lemma~\ref{lem:KKG_C_Ind}.
\end{proof}

\begin{proof}[Proof of Proposition~\ref{pro:mu_forgetful}]
  Since~$X$ is uniformly contractible, it is contractible.  This
  implies that $\RKK_*(X;\C,D)\cong\K_*(D)$.  Using Lemmas
  \ref{identificationofreducedktheoryofbbalgebras}
  and~\ref{homotopyinvarianceofbbbar}, we conclude that the
  $\K$-theory of $\bbbarred(X,D)$ vanishes.  Thus~$\partial'$ is an
  isomorphism.  The maps $\mu_{X,D}^*$ and~$\partial$ are equivalent
  because $X\to\Rips(X)$ is a coarse homotopy equivalence.  We have
  $\partial=\partial'\circ j_*$ by the naturality of the connecting
  map.
\end{proof}

The above description of the co-assembly map in terms of forgetting
control is also available in the equivariant case.  This is a pleasant
and novel feature of our framework.  An analogue of
Proposition~\ref{pro:mu_forgetful} for the coarse Baum-Connes assembly
map is the de\brd{}localisation description in~\cite{Yu}.  Let~$G$ be
a totally disconnected locally compact group with a $G$\nbd{}compact
model for $\EG$.  The space $\EG$ is uniformly contractible and,
moreoever, $H$\nbd{}equivariantly contractible for all compact
subgroups $H\subseteq G$.  Hence, by a small elaboration of the
argument above, $\K_*^H\bigl(\bbbarred_G(\EG,D)\bigr)$ vanishes for
all compact subgroups $H\subseteq G$.  This implies vanishing of
$\Ktop_*\bigl(G,\bbbarred_G(\EG,D)\bigr)$ by~\cite{Ect}.  Hence we
have:

\begin{thm}\label{forgettingcontrol}
  Let~$G$ be a totally disconnected group with a $G$\nbd{}compact
  model for $\EG$.  Then the $G$\nbd{}equivariant coarse
  co\brd{}assembly map for~$G$ is equivalent to the map
  $$
  j_*\colon \Ktop_{*+1}\bigl(G,\bvr_G(\abs\EG,D)\bigr)
  \to\Ktop_{*+1}\bigl(G,\bbred_G(\EG,D)\bigr)
  $$
  induced by the inclusion $j\colon \bvr(\abs\EG,D)
  \to\bbred(\EG,D)$.
\end{thm}

\section{The coarse co-assembly map and equivariant Kasparov theory}
\label{sec:coarse_bivariant}

We first identify the equivariant co-assembly map for~$\abs{G}$ with
coefficients with a map of the form $p_\EG^*\colon
\KK^G(\C,B)\to\RKK^G(\EG;\C,B)$ for suitable~$B$.  Then we prove a
weaker result for general~$B$.

\subsection{An equivalence of maps}
\label{sec:coarse_bivariant_I}

Throughout this section, we fix a locally compact group~$G$, a compact
subgroup $H\subseteq G$, and an $H$\nbd{}$C^*$\brd{}algebra~$D$.  The
induced $G$\nbd{}$C^*$\brd{}algebra $\Ind_H^G(D)$ is defined as
$$
\Ind_H^G D \defeq
\bigl\{f\in C_0(G,D) \bigm|
  \text{$\alpha_h\bigl(f(gh)\bigr)=f(g)$ for all $h\in H$, $g\in G$}
\bigr\},
$$
with $G$\nbd{}action $(gf)(g') = f(g^{-1}g')$.  If~$\Hilm$ is an
$H$\nbd{}equivariant Hilbert $D$\nbd{}module, then a similar formula
defines a $G$\nbd{}equivariant Hilbert $\Ind_H^G(D)$\brd{}module
$\Ind_H^G(\Hilm)$.

We also let $\EG$ be a universal proper $G$\nbd{}space.  Given two
$G$\nbd{}$C^*$\brd{}algebras $A$ and~$B$, we define the bivariant
Kasparov groups $\RKK^G(\EG;A,B)$ as in~\cite{Kasparov} and let
\begin{displaymath}
  p_\EG^*\colon \KK^G(A,B)\to\RKK^G(\EG;A,B)
\end{displaymath}
be the natural map induced by the constant map from $\EG$ to a point.

\begin{thm}  \label{the:equivariant_versus_coarse}
  There are isomorphisms $\KK^G_*(\C,\Ind_H^G D) \cong
  \K_{*+1}^H\bigl(\bvr_H(\abs{G},D)\bigr)$ and $\RKK^G_*(\EG;
  \C,\Ind_H^G D) \cong \KX^*_H(\abs{G},D)$ making the following
  diagram commute:
  $$
  \xymatrix@C=4em{
    \KK^G_*(\C,\Ind_H^G D) \ar[d]^{\cong} \ar[r]^-{p_\EG^*} &
    \RKK^G_*(\EG;\C,\Ind_H^G D) \ar[d]^{\cong} \\
    \K_{*+1}^H\bigl(\bvr_H(\abs{G},D)\bigr)
    \ar[r]^-{\mu^*_{\abs{G},H,D}} &
    \KX^*_H(\abs{G},D).
    }
  $$
\end{thm}

Note that if $D=\C$ with trivial action of~$H$, then $\Ind_H^G D
=C_0(G/H)$.  Hence Theorem~\ref{the:equivariant_versus_coarse}
identifies~\eqref{eq:analytic_coassembly_H} with the
$H$\nbd{}equivariant coarse co\brd{}assembly map $\mu_{\abs{G},H}^*$
and~\eqref{eq:analytic_coassembly} with the coarse co\brd{}assembly
map~$\mu_{\abs{G}}^*$.

We prepare the proof of Theorem~\ref{the:equivariant_versus_coarse}
with several lemmas.  If~$H$ acts on a $C^*$\nbd{}algebra~$A$, let
$A^H$ be the subalgebra of $H$\nbd{}invariant elements of~$A$.

\begin{lemma}  \label{lem:KKG_C_Ind}
  There are natural isomorphisms
  $$
  \KK^G_*(\C,\Ind_H^G D)
  \cong \K_{*+1}(\bvr_H(\abs{G},D)^H)
  \cong \K_{*+1}^H\bigl(\bvr_H(\abs{G},D)\bigr).
  $$
\end{lemma}

\begin{proof}
  We only treat the case $*=0$, the case $*=1$ is similar.  To prove
  the first isomorphism, we describe the cycles for
  $\KK^G_0(\C,\Ind_H^G D)$ more concretely.  Such a cycle is given by
  two $G$\nbd{}equivariant Hilbert modules $\Hilm_\pm$ over $\Ind_H^G
  D$ and a $G$\nbd{}continuous adjointable operator
  $F\colon\Hilm_+\to\Hilm_-$ for which $1-FF^*$, $1-F^*F$ and $gF-F$
  for $g\in G$ are compact.  Let $D_H^\infty\defeq D\otimes
  L^2(H)\otimes\ell^2(\N)$ be the standard $H$\nbd{}equivariant
  Hilbert module over~$D$.  Then $\Ind_H^G(D_H^\infty)$ is naturally
  isomorphic to the standard Hilbert module $\Ind_H^G D\otimes
  L^2(G)\otimes\ell^2(\N)$ over $\Ind_H^G D$.  Since $\Ind_H^G D$ is a
  proper $G$\nbd{}$C^*$\brd{}algebra, the Equivariant Stabilisation
  Theorem of~\cite{Meyer:KKG} applies and yields that every countably
  generated Hilbert module over $\Ind_H^G D$ is absorbed by the
  standard one.  It follows that we can also define
  $\KK^G_0(\C,\Ind_H^G D)$ using only those ``special'' cycles where
  $\Hilm_+=\Hilm_-=\Ind_H^G (D_H^\infty)$.

  Elements of $\Ind_H^G(D_H^\infty)$ are functions in
  $C_0(G,D_H^\infty)$ that satisfy $f(g)=\alpha_h \bigl(f(gh)\bigr)$
  for all $g\in G$, $h\in H$.  The $\Ind_H^G D$\brd{}Hilbert module
  structure is given by pointwise multiplication and pointwise inner
  products.  The group~$G$ acts by left translation.  Thus the space
  of adjointable operators on $\Ind_H^G(D_H^\infty)$ can be identified
  with the space of $*$\nbd{}strictly continuous functions $f\colon
  G\to\Bound(D_H^\infty)$ that are $H$\nbd{}invariant, that is,
  $f(g)=\alpha_h\bigl(f(gh)\bigr)$ for all $g\in G$, $h\in H$.  In
  particular, we can view~$F$ as such a function, which we still
  denote~$F$.

  The $G$\nbd{}continuity of~$F$ means that this function is not just
  strictly continuous: it is uniformly norm continuous, that is,
  \begin{equation}  \label{eq:F_uniformly_continuous}
    \lim_{g\to1} {} \sup_{x\in G} {} \norm{F(g^{-1}x)-F(x)}=0.
  \end{equation}
  Given~\eqref{eq:F_uniformly_continuous}, the condition
  $gF-F\in\Comp(\Ind_H^G D_H^\infty)$ for $g\in G$ translates into the
  two conditions $F(g^{-1}x)-F(x)\in\Comp(D^\infty_H)$ for all $g,x\in
  G$ and
  \begin{equation}  \label{eq:F_vv}
    \lim_{x\to\infty} {} \sup_{g\in K} {} \norm{F(g^{-1}x)-F(x)}=0.
  \end{equation}
  Conversely, \eqref{eq:F_vv} together with ordinary continuity
  implies~\eqref{eq:F_uniformly_continuous}.  Thus we get exactly the
  condition that $F\in\vvr_H(\abs{G},D)^H$.  Since $\Comp(\Ind_H^G
  D^\infty_H) \cong C_0(G,D\otimes\Comp_H)^H$, the compactness of
  $1-FF^*$ and $1-F^*F$ means that the image of~$F$ in
  $\bvr_H(\abs{G},D)^H$ is unitary.  Summing up, ``special'' cycles
  for $\KK^G_0(\C,\Ind_H^G D)$ are in bijection with elements of
  $\vvr_H(\abs{G},D)^H$ whose image in $\bvr_H(\abs{G},D)^H$ is
  unitary.

  Two cycles for $\KK^G_0(\C,\Ind_H^G D)$ differ by a compact
  perturbation if and only if they have the same image in
  $\bvr_H(\abs{G},D)^H$.  The map $\vvr_H(\abs{G},D)^H
  \to\bvr_H(\abs{G},D)^H$ is surjective because~$H$ is compact.
  Therefore, equivalence classes of ``special cycles'' up to compact
  perturbation correspond bijectively to unitaries in
  $\bvr_H(\abs{G},D)^H$.  A cycle is degenerate if and only if it is a
  constant function on~$G$.  Cycles are operator homotopic if and only
  if the resulting unitaries in $\bvr_H(\abs{G},D)^H$ are homotopic.

  It is easy to see that $\bvr_H(\abs{G},D)^H$ is matrix stable.
  Hence we do not have to adjoin matrices to compute its
  $\K$\nbd{}theory.  The subalgebra of constant functions in
  $\bvr_H(\abs{G},D)^H$ is isomorphic to $\Mult(D\otimes\Comp_H)^H$
  and hence has vanishing $\K$\nbd{}theory: the same Eilenberg swindle
  that proves this fact for stable multiplier algebras works
  equivariantly.  As a result, addition of degenerate cycles and
  operator homotopy generate the same equivalence relation on
  ``special'' cycles for $\KK^G_0(\C,\Ind_H^G D)$ as stable homotopy
  equivalence for unitaries in $\bvr_H(\abs{G},D)^H$.  Since operator
  homotopy and homotopy generate the same equivalence relation, we get
  $\KK^G_0(\C,\Ind_H^G D) \cong \K_1(\bvr_H(\abs{G},D)^H)$ as claimed.

  To prove the second isomorphism, we show that $\bvr_H(X,D)^H$ and
  $\bvr_H(X,D)\cross H$ are Morita-Rieffel equivalent (notice that
  both algebras are $\sigma$\nbd{}unital).  If~$H$ is finite, then
  $\Comp(L^2H)$ is finite dimensional, so that
  $$
  \bvr_H(X,D)^H
  \cong \bigl(\bvr_H(X,D)\otimes\Comp(L^2H)\bigr)^H
  \cong \bvr_H(X,D)\cross H
  $$
  by the proof of the Green-Julg Theorem.  Thus we have an
  isomorphism in this case.  For general compact~$H$, we use that the
  fixed point algebra is Morita-Rieffel equivalent to a certain
  ideal~$I$ in the crossed product (see~\cite{Meyer:Fixed}).  The
  imprimitivity bimodule is $\bvr_H(X,D)$ equipped with appropriate
  structure.  We embed $C(H)\subseteq\Bound(L^2H)\subseteq
  \bvr_H(X,D)$ unitally as constant functions on~$X$.  Since this
  embedding is equivariant, the ideal~$I$ contains the corresponding
  ideal for $C(H)$, which is all of $C(H)\cross H\cong\Comp(L^2H)$.
  It follows that~$I$ contains an approximate identity and hence must
  be all of $\bvr_H(X,D)\cross H$.
\end{proof}

We next recall the following well-known facts.

\begin{lemma}  \label{tinylemma}
  Let~$G$ be a locally compact group and~$H$ a closed subgroup.
  Let~$A$ be a $\sigma$\nbd{}$G$\nbd{}$C^*$\brd{}algebra and~$B$ a
  $\sigma$\nbd{}$H$\nbd{}$C^*$\brd{}algebra.
  \begin{enumerate}[\ref{tinylemma}.1.]
  \item The $\sigma$\nbd{}$G$\nbd{}$C^*$\brd{}algebras
    $A\otimes\Ind_H^G B$ and $\Ind_H^G(A\otimes B)$ are isomorphic.

  \item The $\sigma$\nbd{}$C^*$\brd{}algebras $\bigl(\Ind_H^G
    A\bigr)\cross G$ and $A\cross H$ are Morita-Rieffel equivalent.

  \end{enumerate}
\end{lemma}

\begin{cor}  \label{cor:induced_EG}
  Let~$G$ be a discrete group, $H$ a finite subgroup, and~$D$ an
  $H$\nbd{}$C^*$\brd{}algebra.  Then there is a canonical isomorphism
  $$
  \K_*\bigl(C_0(\abs\EG,D)\cross H\bigr) \cong
  \K_*\bigl(C_0(\abs\EG,\Ind_H^G D)\cross G \bigr).
  $$
\end{cor}

\begin{proof}
  Lemma~\ref{tinylemma} implies
  \begin{displaymath}
    C_0(\abs\EG,\Ind_H^G D) \cross G
    \cong \Ind_H^G( C_0(\abs\EG)\otimes D) \cross G
    \sim C_0(\abs\EG,D)\cross H,
  \end{displaymath}
  where~$\cong$ denotes isomorphism and~$\sim$ denotes Morita-Rieffel
  equivalence.
\end{proof}

\begin{lemma}  \label{lem:RKKG_as_Ktheory}
  Let~$G$ be a locally compact group and let~$X$ be a locally compact
  proper $G$\nbd{}space that can be written as a union of an
  increasing sequence $(X_n)$ of $G$\nbd{}compact closed subspaces.
  Let~$A$ be a $C^*$\brd{}algebra with trivial action of~$G$ and
  let~$B$ be a $G$\nbd{}$C^*$\brd{}algebra.  Then there is a natural
  isomorphism
  $$
  \RKK^G_*(X;A,B) \cong \KK_*(A,\varprojlim C_0(X_n,B)\cross G).
  $$
\end{lemma}

Here we use the bivariant Kasparov theory for
$\sigma$\nbd{}$C^*$\brd{}algebras defined by Alexander Bonkat
in~\cite{Bonkat:Thesis}.  We will only apply this lemma for $A=\C$,
where this reduces to $\K$\nbd{}theory for
$\sigma$\nbd{}$C^*$\brd{}algebras as defined in~\cite{Phillips}.

\begin{proof}
  We check that both groups agree on the level of cycles after some
  standard simplifications.  Since $C_0(X,B)$ is a proper
  $G$\nbd{}$C^*$\brd{}algebra, the reduced and full crossed products
  for $C_0(X,B)$ agree.  Moreover, the $C^*$\nbd{}categories of
  $G$\nbd{}equivariant Hilbert modules over $C_0(X,B)$ and of Hilbert
  modules over $C_0(X,B) \cross G$ are equivalent
  (see~\cite{Meyer:Fixed}).  That is, any $G$\nbd{}equivariant Hilbert
  module~$\Hilm$ over $C_0(X,B)$ corresponds to a Hilbert
  module~$\tilde\Hilm$ over $C_0(X,B)\cross G$.  The correspondence is
  such that $\Bound(\tilde\Hilm)$ is naturally isomorphic to the
  $C^*$\nbd{}algebra $\Bound(\Hilm)^G$ of $G$\nbd{}equivariant
  adjointable operators on~$\Hilm$.  The compact operators
  on~$\tilde\Hilm$ correspond to the \emph{generalised fixed point
    algebra} of $\Comp(\Hilm)$, which is the closed linear span of
  operators of the form $\int_G \alpha_g(\ket{\xi}\bra{\eta})\,dg$,
  where $\xi,\eta\in\Hilm$ are \emph{compactly supported} sections.
  (The support of~$\xi$ is the set of $x\in X$ with $\xi_x\neq 0$.)
  More generally, if $T\in\Comp(\Hilm)$ has compact support, then
  $\int_G \alpha_g(T)\,dg$ belongs to the generalised fixed point
  algebra.

  To simplify our notation, we consider the $\sigma$\nbd{}locally
  compact space $\X\defeq \bigcup X_n$ and let $C_0(\X,B)\cross
  G\defeq \varprojlim C_0(X_n,B)\cross G$.  We have natural maps
  $$
  C_0(X,B)\cross G\to C_0(\X,B)\cross G \to C_0(X_n,B)\cross G
  $$
  for all $n\in\N$.  If~$\tilde\Hilm$ is a Hilbert module over
  $C_0(X,B)\cross G$, its restriction to~$X_n$ is the Hilbert module
  $\tilde\Hilm_n\defeq \tilde\Hilm\otimes_{C_0(X,B)\cross G}
  C_0(X_n,B)\cross G$ over $C_0(X_n,B)\cross G$.  Then
  $$
  \hat\Hilm\defeq \varprojlim \tilde\Hilm_n
  \cong \tilde\Hilm \otimes_{C_0(X,B)\cross G} C_0(\X,B)\cross G
  $$
  is a Hilbert module over $C_0(\X,B)\cross G$.  Conversely, given
  a Hilbert module $\hat\Hilm$ over $C_0(\X,B)\cross G$, we obtain a
  Hilbert module over $C_0(X,B)$ by completing the subspace of
  compactly supported sections $\hat\Hilm\cdot C_c(X,B)\cross G$.  It
  is easy to see that these two operations are inverse to each other.
  We have
  $$
  \Comp(\hat\Hilm)\defeq \varprojlim \Comp(\tilde\Hilm_n),
  \qquad
  \Bound(\hat\Hilm)\defeq \varprojlim \Bound(\tilde\Hilm_n).
  $$
  We can describe $\Bound(\tilde\Hilm)$ as the
  $C^*$\nbd{}subalgebra of bounded elements in the
  $\sigma$\nbd{}$C^*$\nbd{}algebra $\Bound(\hat\Hilm)$.  Thus any
  $*$\nbd{}homomorphism $A\to\Bound(\hat\Hilm)$ factorises through
  $\Bound(\tilde\Hilm)$.

  Cycles for $\RKK^G_0(X;A,B)$ are triples $(\Hilm,\phi,F)$
  where~$\Hilm$ is a graded $G$\nbd{}equivariant Hilbert module over
  $C_0(X,B)$, $\phi\colon A\to\Bound(\Hilm)$ is an equivariant
  $*$\nbd{}homomorphism and $F\in\Bound(\Hilm)$ is a self-adjoint,
  $G$\nbd{}equivariant, odd operator for which $M_h\cdot[\phi(a),F]$
  and $M_h\phi(a)(1-F^2)$ are compact for all $h\in C_0(X)$, $a\in A$.
  It is shown in~\cites{Kasparov, Meyer:KKG} that we can arrange~$F$
  to be strictly equivariant, using that~$X$ is proper.  Since~$G$
  acts trivially on~$A$, the range of~$\phi$ consists of
  $G$\nbd{}equivariant operators on~$\Hilm$.

  By our category equivalence, this data is equivalent to a triple
  $(\hat\Hilm,\hat\phi,\hat{F})$, where~$\hat\Hilm$ is a Hilbert
  module over $C_0(\X,B)\cross G$ and $\hat\phi$ and~$\hat{F}$ are
  obtained from $\phi$ and~$F$ using $\Bound(\Hilm)^G
  \cong\Bound(\tilde\Hilm) \subseteq\Bound(\hat\Hilm)$.
  Thus~$\hat\phi$ is a $*$\nbd{}homomorphism and~$\hat{F}$ is a
  bounded, odd, self-adjoint operator.  We claim that this
  construction yields a bijection between cycles $(\Hilm,\phi,F)$ for
  $\RKK^G_0(X;A,B)$ with equivariant~$F$ and cycles
  $(\hat\Hilm,\hat\phi,\hat{F})$ for $\KK_0(A, C_0(\X,B)\cross G)$.
  Let~$S$ be $\phi(a)(F^2-1)$ or $[F,\phi(a)]$ for some $a\in A$ and
  let~$\hat{S}$ be the associated operator on~$\hat\Hilm$.  The proof
  is finished if we show that $M_h S$ is compact for all $h\in C_0(X)$
  if and only if~$\hat{S}$ is compact.

  Assume first that $M_h S\in\Comp(\Hilm)$ for all $h\in C_0(X)$.
  Choose $n\in\N$.  By the properness of the $G$\nbd{}action, there is
  a function $h\in C_c(X)$ with $\int_G h(xg)=1$ for all $x\in X_n$.
  Then~$S$ and $\int_G \alpha_g(M_h S) \,dg$ have the same restriction
  to~$X_n$.  Since $M_h S$ is compact by hypothesis and has compact
  support, this integral belongs to $\Comp(\tilde\Hilm_n)$.  This
  implies $\hat{S}\in\Comp(\hat\Hilm)$.  Suppose conversely that
  $\hat{S}\in\Comp(\hat\Hilm)$ and fix $h\in C_c(X)$.  Choose~$n$ so
  that~$X_n$ contains the support of~$h$.  Thus the product $M_h S$
  only sees the restriction of~$S$ to~$X_n$.  The operator
  on~$\tilde\Hilm_n$ induced by~$\hat{S}$ is compact.  Thus~$S$
  belongs to the generalised fixed point algebra of $\Hilm_n\defeq
  \Hilm\otimes_{C_0(X,B)} C_0(X_n,B)$.  That is, it can be
  approximated by operators of the form $\int_G \alpha_g(T)\,dg$ for a
  finite rank operator~$T$ on~$\Hilm_n$ with compact support.  Hence
  the function $g\mapsto M_h\alpha_g(T)$ has compact support, so that
  $M_h \int_G \alpha_g(T)\,dg$ is a compact operator on~$\Hilm_n$.
  Since these operators approximate $M_h S$, we get $M_h
  S\in\Comp(\Hilm_n)\subseteq\Comp(\Hilm)$.
\end{proof}

\begin{proof}[Proof of Theorem~\ref{the:equivariant_versus_coarse}]
  Lemma~\ref{lem:RKKG_as_Ktheory} and Corollary~\ref{cor:induced_EG}
  yield isomorphisms
  \begin{multline*}
    \RKK^G_*(\EG;\C,\Ind_H^G D)
    \cong
    \K_*(C_0(\abs{\EG},\Ind_H^G D)\cross G)
    \\ \cong
    \K_*(C_0(\abs{\EG},D)\cross H)
    = \KX^*_H(G,D).
  \end{multline*}
  The other isomorphism required for
  Theorem~\ref{the:equivariant_versus_coarse} is provided by
  Lemma~\ref{lem:KKG_C_Ind}.  To check that the resulting diagram
  commutes, we work with (generalised) fixed point algebras instead of
  crossed products, as this simplifies the arguments.  We have
  $C_0(\abs{\EG},D\otimes\Comp_H)^H\cong
  C_0(\abs{\EG},D\otimes\Comp)\cross H$ because $\Comp_H$ contains a
  factor $\Comp(L^2H)$.

  We again define $\KK^G_0(\C,\Ind_H^G D)$ by ``special'' cycles and
  identify them with elements of $F\in\vvr_H(\abs{G},D)^H$ whose image
  $\pi(F)$ in $\bvr_H(\abs{G},D)^H$ is unitary.  Thus the isomorphism
  $\KK^G_0(\C,\Ind_H^G D) \to \K_1(\bvr_H(\abs{G},D)^H)$ maps the
  class represented by the cycle~$F$ to the class represented by the
  unitary $\pi(F)$.

  In order to compute the image of $[\pi(F)]$ under the coarse
  co\brd{}assembly map, we have to describe the connecting map in
  $\K$\nbd{}theory.  Let $0\to I\to E\to Q\to 0$ be an extension of
  $\sigma$\nbd{}$C^*$\brd{}algebras and suppose that $E$ and~$Q$ are
  unital.  Let $u\in Q$ be unitary and lift it to $F\in
  E\subseteq\Mult(I)$.  Since~$u$ is unitary, $FF^*-1$ and $F^*F-1$
  belong to~$I$.  Thus $F\in\Mult(I)$ is a cycle for
  $\KK_0(\C,I)\cong\K_0(I)$.  This element of $\K_0(I)$ is the image
  of $[u]$ under the connecting map.  Thus we get
  $\mu^*_{\abs{G},H,D}[\pi(F)]$ if we lift
  $\pi(F)\in\bvr_H(\abs{G},D)^H\cong\bvr_H(\abs{\EG},D)^H$ to an
  element of $\vvr_H(\abs{\EG},D)^H$ and then view this as a Fredholm
  multiplier of $C_0(\abs{\EG},D\otimes\Comp_H)^H$.

  Choose an $H$\nbd{}invariant continuous function $c\colon
  \EG\to\R_+$ with $\int_G c(xg)\,dg=1$ for all $x\in\EG$ such that
  $S_Y\defeq \supp c\cap Y$ is compact for all $G$\nbd{}compact
  subsets $Y\subseteq\EG$.  Let $L_Y\subseteq G$ be the set of $g\in
  G$ with $S_Y\,g\cap S_Y\neq\emptyset$.  We let
  $$
  \bar{F}(x) \defeq \int_G c(xg)\,F(g^{-1})\,dg.
  $$
  If $x\in S_Y\,g$ for some $g\in G$, then $\bar{F}(x)$ is an
  average of $F(h)$ with $xh^{-1}\in S_Y$, so that $h\in L_Y^{-1}g$.
  It follows that $\bar{F}|_Y$ belongs to $\vvr_H(Y,D)^H$ for all
  $G$\nbd{}compact~$Y$, that is, $\bar{F}\in\vvr_H(\abs\EG,D)^H$.  The
  quotient map $\vvr_H(\abs{\EG},D)^H\to\bvr_H(\abs{G},D)^H$ simply
  restricts a function on~$\EG$ to any $G$\nbd{}orbit in $\abs\EG$.
  Hence $\pi(\bar{F})=\pi(F)$ in $\bvr_H(\abs{G},D)^H$.  Thus
  $\mu^*_{\abs{G},H,D}$ maps $[\pi(F)]$ to the class represented by
  the Fredholm multiplier~$\bar{F}$ of
  $C_0(\abs{\EG},D\otimes\Comp_H)^H$.

  Now go around the diagram the other way.  By definition,
  $p_\EG^*[F]$ is represented by the multiplication operator $F'
  f(x,g)\defeq F(g) f(x,g)$ on $C_0(\EG,\Ind_H^G D_H^\infty)$, with
  action of~$G$ coming from the action on $\EG\times G$ by $h\cdot
  (x,g)=(xh^{-1},hg)$.  The same formulas work if we replace
  $C_0(\EG)$ by $C_0(\abs{\EG})$.  Let~$c$ be as above.  It is easy to
  check that the multiplication operator~$F''$ defined by
  $$
  (F'' f)(x,g)\defeq \int_G c(xh) F(h^{-1}g)\,dh \cdot f(x,g)
  $$
  is a $G$\nbd{}equivariant compact perturbation of~$F'$.  That is,
  $F''$ and~$F'$ have the same class in $\RKK^G_0(\EG;\C,\Ind_H^G D)$
  and~$F''$ is a multiplier of the generalised fixed point algebra of
  $C_0(\abs\EG,\Ind_H^G D_H^\infty)$.  Restriction to
  $\EG\times\{1\}\subseteq\EG\times G$ identifies this generalised
  fixed point algebra with $C_0(\abs\EG,D\otimes\Comp_H)^H$.  The
  isomorphisms
  \begin{multline*}
    \RKK^G(\EG;\C,\Ind_H^G D)\cong
    \K_0(C_0(\abs\EG,\Ind_H^G D)\cross G)
    \\ \cong
    \K_0(C_0(\abs\EG,D\otimes\Comp)\cross H) \cong
    \K_0(C_0(\abs\EG,D\otimes\Comp_H)^H)
  \end{multline*}
  constructed above send~$F'$ to the class of the Fredholm multiplier
  $F''|_{\EG\times\{1\}}$ of $C_0(\abs\EG,D\otimes\Comp_H)^H$.  The
  reason for this is that Lemma~\ref{tinylemma}.2 is proved using the
  same manipulations of generalised fixed point algebras that we used
  above to view~$F''$ as a multiplier of
  $C_0(\abs\EG,D\otimes\Comp_H)^H$.  By construction,
  $F''|_{\EG\times\{1\}}=\bar{F}$.  Thus the diagram commutes as
  desired.
\end{proof}

\subsection{Constructing Kasparov cycles from the stable Higson
  corona}
\label{sec:construct_KK_vv}

\begin{thm}  \label{the:coarse_gives_classes}
  Let $G$ and~$H$ be locally compact groups and let~$X$ be a coarse
  space equipped with commuting actions of $G$ and~$H$.  Suppose
  that~$G$ acts by translations and that~$H$ acts properly and by
  isometries.  Let~$D$ be an $H$\nbd{}$C^*$\brd{}algebra.  Then there
  is a natural commuting diagram
  $$
  \xymatrix@C=4em{
    \K_{*+1}(\bvr_H(X,D)\cross H)
    \ar[r]^-{\mu_{X,H,D}^*} \ar[d]^{\psi^{G,X,H,D}_*} &
    \KX^*_H(X,D) \ar[d]^{\phi^{G,X,H,D}_*} \\
    \KK^G_*(\C,C_0(X,D)\cross H) \ar[r]^-{p_\EG^*} &
    \RKK^G_*(\EG;\C,C_0(X,D)\cross H).
  }
  $$
  If $H\subseteq G$ is a compact subgroup and $X=\abs{G}$ with
  actions of $G$ and~$H$ by multiplication, then this diagram is
  equivalent to the one constructed in
  Theorem~\ref{the:equivariant_versus_coarse} via the Morita-Rieffel
  equivalence $C_0(G,D)\cross H\sim C_0(G,D)^H=\Ind_H^G(D)$.  Hence
  the vertical maps are isomorphisms in this special case.
\end{thm}

We will apply this theorem to construct elements of $\KK^G(\C,B)$ in
Section~\ref{sec:Lipschitz}.  It is also used to construct dual Dirac
morphisms from compactifications
(Theorem~\ref{the:dual_Dirac_compactify}).  Of course, the natural
transformations $\phi$ and~$\psi$ cannot be isomorphisms in general.

\begin{proof}
  We only construct the diagram in the case $*=0$, the case $*=1$ is
  similar.  To simplify notation, let $A\defeq C_0(X,D)\cross H$ and
  assume $D\otimes\Comp_H\cong D$, so that we can omit the
  stabilisations in the definition of $\vvr_H(X,D)$ and $\bvr_H(X,D)$.
  Recall that we have a canonical $C^*$\nbd{}extension
  $$
  0\to A\to \vvr_H(X,D)\cross H \to
  \bvr_H(X,D)\cross H \to 0.
  $$
  Let $u\in\bvr_H(X,D)\cross H$ be unitary and lift it to
  $F\in\vvr_H(X,D)\cross H$.  Then we may view~$F$ as a Fredholm
  multiplier $F\colon A\to A$.  Since $F\in\vvr_H(X,D)$ and~$G$ acts
  on~$X$ by translations, $g\mapsto gF-F$ is a norm continuous
  function from~$G$ to~$A$.  Thus~$F$ defines a cycle for
  $\KK^G_0(\C,A)$.  As in the proof of Lemma~\ref{lem:KKG_C_Ind}, one
  checks that this construction defines a map $\K_1(\bvr_H(X,D)\cross
  H) \to \KK^G_0(\C,A)$.

  There is a function $c\colon \EG\to\R_+$ for which $\int_\EG c(\mu
  g)\,dg=1$ for all $\mu\in\EG$ and $\supp c\cap Y$ is compact for
  $G$\nbd{}compact $Y\subseteq\EG$.  Recall that $\Rips(X)$ contains
  all compactly supported probability measures on~$X$.  Hence we can
  define a map $c_*\colon \EG\times X\to\Rips(X)$ by
  $$
  \langle c_*(\mu,x), \alpha\rangle
  \defeq \int_G c(\mu g) \alpha(g^{-1}x)\,dg
  $$
  for all $\alpha\in C_0(X)$.  One checks easily that this map is
  continuous and satisfies $c_*(\mu g,g^{-1}xh)=c_*(\mu,x)h$ for all
  $g\in G$, $\mu\in\EG$, $x\in X$, $h\in H$.  If $K\subseteq \EG$ is
  compact, then there is a compact subset $L\subseteq G$ such that
  $c(\mu g)=0$ for $\mu\in K$ and $g\notin L$.  Hence $c_*(\mu,x)$ is
  supported in $L^{-1}x$ for $\mu\in S$.  Since~$G$ acts on~$X$ by
  translations, this is contained in $\Subrips_{E,1}(X)$ for a
  suitable entourage~$E$.  The restriction $c_*\colon K\times X\to
  \Subrips_{E,1}(X)$ is proper because it is close to the map
  $$
  K\times X\overset{\pi_X}\longrightarrow
  X\overset{j_X}\longrightarrow \Subrips_{E,1}(X).
  $$

  Recall that $\KX^0_H(X,D)\defeq \K_0(C_0(\Rips(X),D)\cross H)$.  We
  briefly write~$B$ for the $\sigma$\nbd{}$C^*$\brd{}algebra
  $C_0(\Rips(X),D)\cross H$.  Elements of $\K_0(B)$ are represented by
  (bounded) Fredholm multipliers $F\in\Mult(B)$ (that is, $1-F^*F\in
  B$ and $1-FF^*\in B$).  We do not have to stabilise because~$B$ is
  already stable.  View such a Fredholm multiplier as a function
  $H\times \Rips(X)\to\Mult(D)$ and pull it back with $c_*\colon
  \EG\times X\to\Rips(X)$ to a function $c^*(F)\colon H\times
  \EG\times X\to\Mult(D)$.  The equivariance of~$c_*$ implies that
  $c^*(F)$ is a $G$\nbd{}invariant multiplier of $C_0(\EG,A)$.  The
  restriction of $c^*(F)$ to $C(K,A)$ is Fredholm for all compact
  subsets $K\subseteq\EG$ because~$c_*$ restricts to a proper map
  $K\times X\to\Subrips_{E,1}(X)$ for some entourage~$E$.  Therefore,
  $c^*(F)$ is a cycle for $\RKK^G_0(\EG;\C,A)$.  This yields a natural
  map $\KX^0_H(X,D)\to\RKK^G_0(\EG;\C,A)$.

  Finally, routine computations which we leave to the reader show that
  the diagram in the statement of the theorem commutes and agrees with
  the one in Theorem~\ref{the:equivariant_versus_coarse} for
  $X=\abs{G}$.
\end{proof}

For any pair of $C^*$\nbd{}algebras $B,P$ we have a canonical exact
sequence
$$
0 \to (C_0(\abs{\EG},B)\otimes P) \cross G
\to (\vvr(\abs{\EG},B)\otimes P) \cross G
\to (\bvr(\abs{G},B) \otimes P) \cross G \to 0.
$$
Here~$\otimes$ denotes maximal tensor products.  Let
$\partial_{\abs{G},B,P}$ denote the corresponding boundary map.

\begin{cor}  \label{cor:nat_trafo}
  There is a natural transformation
  $$
  \Psi^{B,P}_*\colon
  \K_{*+1}\bigl((\bvr(\abs{G},B)\otimes P)\cross G\bigr) \to
  \KK_*^G(\C,B\otimes P)
  $$
  for pairs of $G$\nbd{}$C^*$\brd{}algebras $B$, $P$, which makes the
  following diagram commute:
  $$
  \xymatrix@C=4em{
    \K_{*+1}\bigl((\bvr_G(\abs{G},B)\otimes P)\cross G\bigr)
    \ar[r]^-{\partial_{\abs{G},B,P}}
    \ar[d]^{\Psi^{B,P}_*} &
    \K_*(C_0(\abs\EG,B\otimes P)\cross G) \ar[d]_{\cong} \\
    \KK_*^G(\C,B\otimes P) \ar[r]^-{p_\EG^*} &
    \RKK_*^G(\EG;\C,B\otimes P).
  }
  $$
  The map $\Psi^{B,P}_*$ is an isomorphism if $P=C_0(G/H)$ for a
  compact open subgroup $H\subseteq G$.
\end{cor}

We shall see later that $\Psi^{B,\ADir}_*$ is an isomorphism for the
source of the Dirac morphism~$\ADir$ provided there is a
$G$\nbd{}compact model for $\EG$ and~$G$ is almost totally
disconnected (Theorem~\ref{the:Kreg}).

\begin{proof}
  Suppose first that $P=\C$.  Let~$G$ act on~$\abs{G}$ on both sides
  by multiplication and identify $C_0(\abs{G},B)\cross G\cong
  B\otimes\Comp(L^2G)\sim B$ as usual.  The map
  $\psi^{G,\abs{G},G,B}_*$ of Theorem~\ref{the:coarse_gives_classes}
  is the required map for $P=\C$.  One checks easily that the map
  $$
  \phi^{G,\abs{G},G,B}_*\colon \K_*(C_0(\abs\EG,B)\cross G)
  \to \RKK_*^G(\EG;\C,B)
  $$
  is the canonical isomorphism that we have already used above.  To
  get the map for arbitrary~$P$, simply use the embedding
  $\bvr_G(\abs{G},B)\otimes P\to\bvr_G(\abs{G},B\otimes P)$.  If
  $P=C_0(G/H)$ for a compact subgroup, then $A\otimes P\cross G\sim
  A\cross H$ for all $G$\nbd{}$C^*$\brd{}algebras~$A$.  If $H\subseteq
  G$ is open, then there is no difference between $H$\nbd{}continuity
  and $G$\nbd{}continuity.  Hence the second assertion follows from
  the last assertion of Theorem~\ref{the:coarse_gives_classes}.
\end{proof}

\section{Projective resolutions, Dirac and dual Dirac morphisms}
\label{sec:resolutions_Dirac_dual_Dirac}

We recall some results from~\cite{MeyerNest} concerning Dirac and dual
Dirac morphisms and the Baum-Connes assembly map.

Let~$G$ be a locally compact group and $H$ a compact subgroup of~$G$.
We have the \emph{restriction functor} $\Res_G^H\colon \KK^G\to\KK^H$,
whose definition is obvious, and the \emph{induction functor}
$\Ind_H^G\colon\KK^H\to\KK^G$, which we have already used above.  We
call $G$\nbd{}$C^*$\brd{}algebras of the form $\Ind_H^G D$ for
compact~$H$ \emph{compactly induced}.  Let $\CI\subseteq\KK^G$ be the
class of compactly induced $G$\nbd{}$C^*$\brd{}algebras and let
$\gen{\CI}\subseteq\KK^G$ be the localising \label{hello!}subcategory
generated by $\CI$.  This is the smallest full subcategory of $\KK^G$
containing $\CI$ that satisfies
\begin{enumerate}[(1)]
\item $\gen{\CI}$ is triangulated, that is, closed under suspensions
  and under extensions with an equivariant, completely positive,
  contractive section;

\item $\gen{\CI}$ is closed under countable direct sums.

\end{enumerate}
All proper $G$\nbd{}$C^*$\brd{}algebras belong to $\gen{\CI}$ by
\cite{MeyerNest}*{Corollary 7.3}.

An element $f\in\KK^G(A,B)$ is called a \emph{weak equivalence} if
$\Res_G^H (f)$ is invertible in $\KK^H(A,B)$ for all compact subgroups
$H\subseteq G$.  An object $A\in\KK^G$ is called \emph{weakly
  contractible} if $\Res_G^H(A)\cong0$ in $\KK^H$.

\begin{defn}[see \cite{MeyerNest}*{Definition 4.5}]  \label{def:Dirac}
  A \emph{Dirac morphism} for a locally compact group~$G$ is a weak
  equivalence $\Dirac\in\KK^G(\ADir,\C)$ with $\ADir\in\gen{\CI}$.
\end{defn}

Any group~$G$ has a Dirac morphism (\cite{MeyerNest}*{Proposition
  4.6}).  It is unique in the sense that if $\Dirac\in\KK^G(\ADir,\C)$
and $\Dirac'\in\KK^G(\ADir',\C)$ are Dirac morphisms, then there is an
isomorphism $i\colon \ADir\to\ADir'$ with $\Dirac'\circ i=\Dirac$.
From now on we fix a Dirac morphism $\Dirac\in\KK^G(\ADir,\C)$.  For
any $A\in\KK^G$, we have $A\otimes\ADir\in\gen{\CI}$, and
$\ID_A\otimes\Dirac\in\KK^G(A\otimes\ADir,A)$ is a weak equivalence.
Thus $A\otimes\ADir$ is a $\gen{\CI}$\brd{}simplicial approximation
of~$A$.  The morphism $\ID_A\otimes\Dirac$ is invertible if and only
if $A\in\gen{\CI}$.

\begin{thm}[see \cite{MeyerNest}*{Theorem 5.2, Proposition 10.2}]
  \label{the:BC_assembly}
  The Baum-Connes assembly map with coefficients~$A$ is equivalent to
  the map
  $$
  \Dirac_*\colon \K_*((\ADir\otimes A)\rcross G)\to\K_*(A\rcross G),
  $$
  induced by a Dirac morphism $\Dirac\in\KK^G(\ADir,\C)$.
  
  Moreover, the natural projection induces an isomorphism
  $$
  \K_*((\ADir\otimes A)\cross G)
  \cong \K_*((\ADir\otimes A)\rcross G).
  $$
\end{thm}

As a consequence, the functor $A\mapsto\Ktop(G,A)$ is the
\emph{localisation} (or left derived functor) of both $A\mapsto
\K(A\rcross G)$ and $A\mapsto \K(A\cross G)$.  This justifies calling
the map $\Dirac_*\colon F(A\otimes\ADir)\to F(A)$ for a covariant
functor~$F$ the \emph{assembly map for~$F$}.  For a contravariant
functor~$F$, we obtain a \emph{co-assembly map} $F(A)\to
F(A\otimes\ADir)$.

We shall be particularly interested in the contravariant functor
$A\mapsto \KK^G(A,B)$ for $B\in\KK^G$.  Its localisation $A\mapsto
\KK^G(A\otimes\ADir,B)$ also gives the morphisms $A\to B$ in the
localisation of the category $\KK^G$ at the weak equivalences.  The
assembly map for this functor can be described in more classical terms
as follows:

\begin{thm}[see \cite{MeyerNest}*{Theorem 7.1}]
  \label{the:derived_category}
  Let $\EG$ be a locally compact model for the universal proper
  $G$\nbd{}space.  There is a natural isomorphism $\KK^G(\ADir\otimes
  A,B) \cong \RKK^G(\EG;A,B)$ making the following diagram commute:
  \begin{displaymath}
    \xymatrix@C=-1em{
      {\KK^G(A\otimes\ADir,B)} \ar[rr]^-{\cong} & &
      {\RKK^G(\EG;A,B)} \\
      & {\KK^G(A,B).} \ar[ur]_{p_\EG^*}
      \ar[ul]^{\Dirac^*} &
    }
  \end{displaymath}
\end{thm}

It is useful to examine this in greater detail, considering the
diagrams
\begin{gather}
  \label{eq:KK_diagram}
  \begin{gathered}
  \xymatrix{
    \KK^G(A,B) \ar[r]^{\Dirac^*} \ar[dr]^{\tau_\ADir} &
    \KK^G(A\otimes\ADir,B) \\
    \KK^G(A,B\otimes\ADir) \ar[u]^{\Dirac_*} \ar[r]^{\Dirac^*} &
    \KK^G(A\otimes\ADir,B\otimes\ADir), \ar[u]^{\Dirac_*}_{\cong}
  }
  \end{gathered}
  \\
  \label{eq:RKK_diagram}
  \begin{gathered}
  \xymatrix{
    \RKK^G(\EG;A,B) \ar[r]^{\Dirac^*}_{\cong}
    \ar[dr]^{\tau_\ADir}_{\cong} &
    \RKK^G(\EG;A\otimes\ADir,B) \\
    \RKK^G(\EG;A,B\otimes\ADir) \ar[u]^{\Dirac_*}_{\cong}
    \ar[r]^{\Dirac^*}_{\cong} &
    \RKK^G(\EG;A\otimes\ADir,B\otimes\ADir). \ar[u]^{\Dirac_*}_{\cong}
  }
  \end{gathered}
\end{gather}
Here $\tau_\ADir$ denotes the exterior product with~$\ADir$.  In
addition, the maps $p_\EG^*$ give a natural transformation between the
two diagrams.  That is, we get a commuting diagram in the form of a
cube, which we do not draw.

The map $\KK^G(A\otimes\ADir,B\otimes\ADir)\to\KK^G(A\otimes\ADir,B)$
in~\eqref{eq:KK_diagram} is an isomorphism by
\cite{MeyerNest}*{Proposition 4.4}.  The composition
$\Dirac_*\circ\tau_\ADir$ in~\eqref{eq:KK_diagram} agrees
with~$\Dirac^*$ by well-known properties of the exterior product.
Since~$\Dirac_*$ is an isomorphism and the square
in~\eqref{eq:KK_diagram} evidently commutes, the lower triangle also
commutes.  The same argument shows that~\eqref{eq:RKK_diagram}
commutes.  All the maps in~\eqref{eq:RKK_diagram} are isomorphisms
because $p_\EG^*(\Dirac)$ is invertible by \cite{MeyerNest}*{Corollary
7.3}.  Theorem~\ref{the:derived_category} together with invertibility
of $p_\EG^*(\Dirac)$ implies that the map $p_\EG^*$ is an isomorphism
from $\KK^G(A\otimes\ADir,B)$ to $\RKK^G(\EG;A\otimes\ADir,B)$ and,
similarly, from $\KK^G(A\otimes\ADir,B\otimes\ADir)$ to
$\RKK^G(\EG;A\otimes\ADir,B\otimes\ADir)$.

Finally, we observe that the maps in~\eqref{eq:KK_diagram} are all
isomorphisms for $A\in\gen{\CI}$ because then
$\ID_A\otimes\Dirac\in\KK^G(A\otimes\ADir,A)$ is invertible.  For the
same reason, the vertical maps in~\eqref{eq:KK_diagram} are
isomorphisms if $B\in\gen{\CI}$.

\begin{defn}[see \cite{MeyerNest}*{Definition 8.1}]
  \label{def:dual_Dirac}
  If $\eta\in\KK^G(\C,\ADir)$ satisfies $\eta\circ\Dirac = 1_\ADir$,
  then we call~$\eta$ a \emph{dual Dirac morphism} and $\gamma =
  \Dirac\circ\eta$ a \emph{$\gamma$\nbd{}element} for~$G$.
\end{defn}

If a dual Dirac morphism exists, then it is unique.  It is shown in
\cite{MeyerNest}*{Theorem 8.2} that a dual Dirac morphism exists
whenever the group has a $\gamma$\nbd{}element according to one of the
more traditional definitions (see \cite{KasparovSkandalis2}).
Conversely, if~$\ADir$ is a proper $G$\nbd{}$C^*$\brd{}algebra, then
Definition~\ref{def:dual_Dirac} is equivalent to the traditional ones.
We show in Section~\ref{sec:detect_Dirac} that~$\ADir$ can be taken to
be a proper $G$\nbd{}$C^*$\brd{}algebra for many groups; in
particular, this holds for discrete groups with finite dimensional
$\EG$.

It follows from \cite{MeyerNest}*{Theorem 8.3.5} and the above
discussion that the following assertions are all equivalent to the
existence of a dual Dirac morphism:
\begin{enumerate}[(1)]
\item the map $p_\EG^*\colon \KK^G(A,B)\to\RKK^G(\EG;A,B)$ is an
  isomorphism for all $A\in\KK^G$, $B\in\gen{\CI}$;
  
\item the map $\Dirac^*\colon \KK^G(A,B) \to\KK^G(A\otimes\ADir,B)$
  in~\eqref{eq:KK_diagram} is an isomorphism for all $A\in\KK^G$,
  $B\in\gen{\CI}$;
  
\item the map $\Dirac^*\colon \KK^G(A,B\otimes\ADir)
  \to\KK^G(A\otimes\ADir,B\otimes\ADir)$ in~\eqref{eq:KK_diagram} is
  an isomorphism for all $A,B\in\KK^G$.

\end{enumerate}

Using the characterisation $(1)$ and the identification of the coarse
co\brd{}assembly map in Theorem~\ref{the:equivariant_versus_coarse},
we obtain the following theorem:

\begin{thm}  \label{the:dual_Dirac_gives_isomorphism}
  If~$G$ has a dual Dirac morphism (or merely an approximate dual
  Dirac morphism), then the equivariant coarse co\brd{}assembly map
  with coefficients
  $$
  \mu_{\abs{G},H,D}^*\colon
  \K_{*+1}^H\bigl(\bvr_H(\abs{G},D)\bigr) \to
  \KX^*_H(\abs{G},D)
  $$
  is an isomorphism for all compact subgroups $H\subseteq G$ and
  all $H$\nbd{}$C^*$\brd{}algebras~$D$.  In particular,
  $\mu_{\abs{G}}^*\colon \K_{*+1}\bigl(\bvr(\abs{G})\bigr)\to
  \KX^*(\abs{G})$ is an isomorphism.
\end{thm}

We do not define approximate dual Dirac morphisms here; see
\cite{MeyerNest}*{Section 8.1} for a discussion.  A group that acts
properly by isometries on a weakly bolic, weakly geodesic space has an
approximate dual Dirac morphism by a result of Gennadi
Kasparov and Georges Skandalis (\cite{KasparovSkandalis2}).  Thus
Theorem~\ref{the:dual_Dirac_gives_isomorphism} implies that the coarse
co\brd{}assembly map is an isomorphism for such groups.

In the following, we investigate whether the converse of
Theorem~\ref{the:dual_Dirac_gives_isomorphism} holds, that is, whether
isomorphism of $\mu_{\abs{G},H,D}$ for all $H$, $D$ implies the
existence of a dual Dirac morphism for~$G$.  We will use the following
lemma.

\begin{lemma}  \label{lem:dual_Dirac_isomorphism}
  A dual Dirac morphism for~$G$ exists if and only if the map
  \begin{equation}  \label{eq:dual_Dirac_iso}
    p_\EG^*\colon \KK^G(\C,\ADir)\to\RKK^G(\EG;\C,\ADir)
  \end{equation}
  is surjective.  If this is the case, then~\eqref{eq:dual_Dirac_iso}
  is an isomorphism.
\end{lemma}

\begin{proof}
  Theorem~\ref{the:derived_category}
  identifies~\eqref{eq:dual_Dirac_iso} with $\Dirac^*\colon
  \KK^G(\C,\ADir)\to\KK^G(\ADir,\ADir)$.  This yields the first
  assertion because any pre-image of $1_\ADir$ is a dual Dirac
  morphism.  Conversely, if~$G$ has a dual Dirac morphism
  then~\eqref{eq:dual_Dirac_iso} is an isomorphism by
  characterisation~(1) for existence of a dual Dirac morphism above.
\end{proof}

Since the category $\gen{\CI}$ to which~$\ADir$ belongs is generated
by $\CI$, one might hope that isomorphism of $\mu_{\abs{G},H,D}$ for
all $H$, $D$ implies that~\eqref{eq:dual_Dirac_iso} is an isomorphism,
so that~$G$ has a dual Dirac morphism by
Lemma~\ref{lem:dual_Dirac_isomorphism}.  Unfortunately, we know
nothing about the behaviour of $B\mapsto\KK^G(\C,B)$ for infinite
direct sums, and these are needed to construct general objects of
$\gen{\CI}$.  In order to construct~$\ADir$ without using infinite
direct sums, we need finiteness hypotheses on $\EG$.

\section{More about the Dirac dual Dirac method}
\label{sec:more_dual_Dirac}

This section contains some new results about the Dirac dual Dirac
method.  We first explain how to construct concrete models for Dirac
morphisms.  Then we prove that the existence of a dual Dirac morphism
is hereditary for extensions.  It follows that a locally compact group
contains a dual Dirac morphism if and only if its group of connected
components has one.

\subsection{Detecting Dirac morphisms}
\label{sec:detect_Dirac}

\begin{lemma}  \label{lem:poincareduality}
  Let~$G$ be a locally compact group, let~$A$ be a
  $G$\nbd{}$C^*$\brd{}algebra and let $d\in\KK^G(A,\C)$.  Then~$d$ is
  a Dirac morphism for~$G$ if and only if there are natural
  isomorphisms $\KK^G(A,B) \cong \RKK^G(\EG;\C,B)$ for all
  $G$\nbd{}$C^*$\brd{}algebra~$B$ that make the following diagram
  commute:
  \begin{displaymath}
    \xymatrix@C=-1em{
      {\KK^G(A,B)} \ar[rr]^-{\cong} & &
      {\RKK^G(\EG;\C,B)} \\
      & {\KK^G(\C,B).}  \ar[ur]_{p_\EG^*} \ar[ul]^{d^*} &
    }
  \end{displaymath}
\end{lemma}

\begin{proof}
  Theorem~\ref{the:derived_category} shows that a Dirac morphism
  $\Dirac\in\KK^G(\ADir,\C)$ has these properties.  Conversely, the
  hypotheses on~$d$ determine the functor $B\mapsto \KK^G(A,B)$ and
  the natural transformation $d^*\colon \KK^G(\C,B)\to\KK^G(A,B)$
  uniquely.  By the Yoneda Lemma, this implies that $d$ and~$\Dirac$
  are equivalent.
\end{proof}

\begin{cor}  \label{cor:Dirac_EG_manifold}
  Suppose that~$\EG$ can be realised by a proper isometric action
  of~$G$ on a complete Riemannian manifold~$M$.  Then the class
  $[\Dirac_M] \in\KK^G(C_\tau(M),\C)$ constructed by Gennadi Kasparov
  in~\cite{Kasparov}*{Definition 4.2} is a Dirac morphism for~$G$.
\end{cor}

\begin{cor}  \label{cor:Dirac_EG_simplicial}
  Suppose that~$\EG$ can be realised by a finite dimensional
  simplicial complex~$X$ on which~$G$ acts simplicially.  Then the
  class $[\Dirac_X]\in \KK^G(\mathcal{A}_X,\C)$ constructed by Gennadi
  Kasparov and Georges Skandalis in
  \cite{KasparovSkandalis}*{Definition 1.3} is a Dirac morphism
  for~$G$.
\end{cor}

\begin{proof}
  The sufficient condition of Lemma~\ref{lem:poincareduality} is
  verified in \cite{Kasparov}*{Theorem 4.9} and
  \cite{KasparovSkandalis}*{Theorem 6.5}.
\end{proof}

\begin{rmk}
  Formally, the source~$\ADir$ of the Dirac morphism cannot be graded
  because the Kasparov category of graded $C^*$\brd{}algebras is not
  triangulated.  However, it is permissible to use a graded
  $G$\nbd{}$C^*$\brd{}algebra that is $\KK^G$\brd{}equivalent to an
  ungraded one.  It is well-known that $C_\tau(M)$ in
  Corollary~\ref{cor:Dirac_EG_manifold} is $\KK^G$\brd{}equivalent to
  $C_0(T^*M)$.  A similar ungraded model for $\mathcal{A}_X$ is
  constructed in~\cite{KasparovSkandalis}.  Therefore, we may ignore
  this technical issue.
\end{rmk}

We call a locally compact group~$G$ \emph{almost totally disconnected}
if the connected component of the identity element in~$G$ is compact.
Of course, totally disconnected groups have this property.  A group is
almost totally disconnected if and only if there exists a proper
simplicial action of~$G$ on a simplicial complex.  In this case, we
can always realise $\EG$ by a simplicial action on a simplicial
complex.  However, $\EG$ need not be finite dimensional.  This is the
only obstruction to applying Corollary~\ref{cor:Dirac_EG_simplicial}.

\subsection{Dual Dirac morphisms for group extensions}
\label{sec:permanence_dual_Dirac}

Let $N\into E\prto G$ be an extension of locally compact groups.  A
subgroup $U\subseteq E$ is called \emph{$N$\nbd{}compact} if its image
in~$G$ is compact.  Then~$U$ is an extension of~$N$ by a compact
group.

\begin{thm}  \label{the:dual_Dirac_extension}
  Suppose that~$G$ and all $N$\nbd{}compact subgroups of~$E$ have dual
  Dirac morphisms.  Then~$E$ has a dual Dirac morphism as well.
\end{thm}

\begin{proof}
  We assume that~$G$ is almost totally disconnected for simplicity.
  This special case implies the general assertion using
  Corollary~\ref{cor:dual_Dirac_td} (which only requires the special
  case).  Let $\Dirac_G\in\KK^G(\ADir_G,\C)$ and
  $\Dirac_E\in\KK^E(\ADir_E,\C)$ be Dirac morphisms for $G$ and~$E$
  and let $\eta_G\in\KK^G(\C,\ADir_G)$ be a dual Dirac morphism
  for~$G$.

  The homomorphism $\pi\colon E\to G$ induces a functor $\pi^*\colon
  \KK^G\to\KK^E$ satisfying $\pi^*(\C)=\C$.  The functor~$\pi$ maps
  weak equivalences to weak equivalences, since it maps compact
  subgroups to compact subgroups.  Since weak equivalences between
  objects of $\gen{\CI}$ are invertible, $\ID\otimes\pi^*(\Dirac_G)\in
  \KK^E(\ADir_E\otimes\pi^*(\ADir_G), \ADir_E)$ is invertible.  We
  claim that
  \begin{equation}  \label{eq:dual_Dirac_extension}
    (\Dirac_E\otimes\ID)^*\colon
    \KK^E(\pi^*(A),\ADir_E)
    \to \KK^E(\ADir_E\otimes \pi^*(A), \ADir_E)
  \end{equation}
  is an isomorphism for $A=\ADir_G$.  Before we prove this claim, we
  show how it yields a dual Dirac morphism for~$E$.  Let
  $\beta\in\KK^E(\pi^*(\ADir_G),\ADir_E)$ be the pre-image of
  $\ID\otimes\pi^*(\Dirac_G)$ under the
  isomorphism~\eqref{eq:dual_Dirac_extension} and let
  $\eta_E\defeq\beta\circ\pi^*(\eta_G)\in\KK^G(\C,\ADir_E)$.  Then
  $\eta_E\circ\Dirac_E=\ID_{\ADir_E}\otimes\pi^*(\Dirac_G\eta_G)$.
  Hence $\eta_E\circ\Dirac_E$ is an idempotent weak equivalence
  $\ADir_E\to\ADir_E$ and so must be equal to~$1$.  Consequently,
  $\eta_E$ is a dual Dirac morphism for~$E$.

  The class of~$A$ for which~\eqref{eq:dual_Dirac_extension} is an
  isomorphism is localising (see page~\pageref{hello!}).  Therefore,
  it suffices to prove that~\eqref{eq:dual_Dirac_extension} is an
  isomorphism for compactly induced~$A$.  Equivalently, we need only
  consider the case where $\pi^*(A)=\Ind_U^E(D)$ for some
  $N$\nbd{}compact subgroup $U\subseteq E$.  Any compact subgroup
  of~$G$ is contained in a compact open one because~$G$ is almost
  totally disconnected.  Making~$U$ larger, we may assume that
  $\pi(U)\subseteq G$ is open and compact.

  We identify $\ADir_E\otimes\Ind_U^E(D) \cong \Ind_U^E(\ADir_E\otimes
  D)\cong \Ind_U^E(\ADir_U\otimes D)$ because the restriction of a
  Dirac morphism for~$E$ is a Dirac morphism for~$U$ by
  \cite{MeyerNest}*{Proposition 10.1}.  We rewrite $\KK^E(\Ind_U^E
  A,B) \cong \KK^U(A,\Res_E^U B)$ as in \cite{MeyerNest}*{Proposition
  3.1}.  We have $\Res_E^U \ADir_E\in\gen{\CI}$ by
  \cite{MeyerNest}*{Proposition 10.1}.  Since~$U$ has a dual Dirac
  morphism by hypothesis, the map
  $$
  \Dirac_U^*\colon \KK^U(D,\Res_E^U \ADir_E)
  \to \KK^U(\ADir_U\otimes D,\Res_E^U \ADir_E)
  $$
  is an isomorphism by \cite{MeyerNest}*{Theorem 8.3}.  This means
  that~\eqref{eq:dual_Dirac_extension} is an isomorphism for compactly
  induced~$A$ and hence for $A=\ADir_G$.
\end{proof}

\begin{cor}  \label{cor:dual_Dirac_td}
  Let~$G$ be a locally compact group, let $G_0\subseteq G$ be the
  connected component of the identity, and let $G/G_0$ be its group of
  connected components.  Then~$G$ has a dual Dirac morphism if and
  only if $G/G_0$ has one.
\end{cor}

Thus when constructing dual Dirac morphisms we may always restrict
attention to totally disconnected groups.

\begin{proof}
  All $G_0$\nbd{}compact subgroups of~$G$ are almost connected and
  hence have a dual Dirac morphism by~\cite{Kasparov}.  The assertion
  now follows from Theorem~\ref{the:dual_Dirac_extension}.
\end{proof}

%% \begin{rmk}
%%   Let $\EG$ be a universal proper $G$\nbd{}space.  One can construct a
%%   natural isomorphism
%%   $$
%%   \KK^E(A\otimes\pi^*(\ADir_G),B)\cong \RKK^E(\EG;A,B).
%%   $$
%%   as in the proof of \cite{MeyerNest}*{Theorem 7.1}.  If we only
%%   assume the existence of dual Dirac morphisms for $N$\nbd{}compact
%%   subgroups in~$E$, we can still prove by the same method as above
%%   that
%%   $$
%%   \RKK^E(\EG;A,B)\to\RKK^E(\mathcal{E}E;A,B)
%%   $$
%%   is an isomorphism for all $B\in\gen{\CI}$, $A\in\KK^E$.  We get
%%   $$
%%   \RKK^E(\C;A,B)\to\RKK^E(\EG;A,B)
%%   $$
%%   for all $B\in\gen{\CI}$, $A\in\KK^E$ if there is a dual Dirac
%%   morphism for~$G$.
%% \end{rmk}

\section{Descent}
\label{sec:descent}

In this section, we state and prove our Descent Principle.
  Fix a locally compact group~$G$ and a Dirac morphism
$\Dirac\in\KK^G(\ADir,\C)$ for~$G$.  We first treat the case where~$G$
has a $G$\nbd{}compact model for $\EG$.  Then we treat the case
where~$G$ is discrete, does not have too many finite subgroups, and
has a finite dimensional model for $\EG$.

In both cases, we need some information on~$\ADir$, the domain of the
Dirac morphism.  In the case of $G$\nbd{}compact $\EG$, we obtain this
information using generalities on compactly generated triangulated
categories.  For discrete~$G$ with finite dimensional $\EG$, we use
instead the concrete description of~$\ADir$ in
Corollary~\ref{cor:Dirac_EG_simplicial}.  The idea is that we want to
use bootstrapping arguments to obtain assertions about~$\ADir$ from
assertions about much simpler coefficient algebras.

Let~$\GEN$ be some class of objects in $\KK^G$.  Write $\gen{\GEN}$
for the localising subcategory and $\gen{\GEN}_\fin$ for the thick
triangulated subcategory generated by it.  The latter is the smallest
subcategory of $\KK^G$ containing~$\GEN$ that is closed under
suspensions, admissible extensions and retracts.  The larger category
$\gen{\GEN}$ has the same properties and is also closed under
countable direct sums.  In our applications, we know that something
desirable happens for objects of~$\GEN$, and therefore that this also
happens for objects of $\gen{\GEN}_\fin$ for purely formal reasons.
We would therefore like~$\ADir$, which is our target, but which
\emph{a priori} only lies in $\gen{\GEN}$, to actually lie in
$\gen{\GEN}_\fin$.  For this we need a hypothesis on $\EG$.

\subsection{The case of finite classifying space}
\label{sec:finite_EG}

The following is our Descent Principle for groups~$G$ admitting a
$G$\nbd{}compact model for $\EG$.

\begin{thm}  \label{the:dual_Dirac_finite}
  Let~$G$ be a locally compact group with $G$-compact $\EG$.  Then~$G$
  has a dual Dirac morphism if and only if the $H$\nbd{}equivariant
  coarse co\brd{}assembly map
  $$
  \mu_{\abs{G},H}^*\colon \K_{*+1}(\bvr_H(\abs{G})\cross H)
  \to \KX^*_H(\abs{G})
  $$
  is an isomorphism for all smooth compact subgroups $H\subseteq G$.
\end{thm}

\begin{cor}
  If~$G$ is a torsion free discrete group with finite classifying
  space $BG$, then~$G$ has a dual Dirac morphism if and only if the
  coarse co\brd{}assembly map
  $$
  \mu_{\abs{G}}^*\colon \K_{*+1}\bigl(\bvr(\abs{G})\bigr)
  \to \KX^*(\abs{G})
  $$
  is an isomorphism.  In particular, the existence of a dual Dirac
  morphism for~$G$ is a coarse property, that is, it only depends on
  the coarse space~$\abs{G}$.
\end{cor}

The proof of Theorem~\ref{the:dual_Dirac_finite} requires some
preparation.  Assume first that~$G$ is almost totally disconnected.
In this case, any compact subgroup is contained in a compact open
subgroup.  We let
$$
\CI_0\defeq
\{ C_0(G/H) \mid \text{$H\subseteq G$ compact open subgroup}\}.
$$
If~$G$ is a arbitrary, we use instead the larger class of smooth
compact subgroups.  The following definition is equivalent to the one
in~\cite{MeyerNest}.  We call a compact subgroup $H\subseteq G$
\emph{smooth} if there are an open almost connected subgroup
$U\subseteq G$ and a compact normal subgroup~$N$ in~$U$ such that
$U/N$ is a Lie group and $N\subseteq H\subseteq U$.  Then $G/H$ is a
smooth manifold, being a disjoint union of copies of the homogeneous
space $(U/N)/(H/N)$.  We let
$$
\CI_1\defeq
\{ C_0(G/H) \mid \text{$H\subseteq G$ smooth compact subgroup}\}.
$$
If~$G$ is almost totally disconnected, then $\CI_0\subseteq\CI_1$.
If~$G$ is totally disconnected, then $\CI_0=\CI_1$.

We have $\ADir\in\gen{\CI_1}$ for all~$G$ by
\cite{MeyerNest}*{Proposition 9.2}.  If~$G$ is almost
totally disconnected, then $\ADir\in\gen{\CI_0}$.  We are going to
find a criterion for~$\ADir$ to belong to $\gen{\CI_1}_\fin$ or
$\gen{\CI_0}_\fin$ that uses abstract results on triangulated
categories.  These are related to the Brown Representability Theorem
with cardinality restrictions \cite{MeyerNest}*{Theorem 6.1}.

We call $A\in\KK^G$ \emph{compact} if $\KK^G(A,B)$ is countable for
all $B\in\KK^G$ and, in addition, $B\mapsto\KK^G(A,B)$ commutes with
countable direct sums.

\begin{lemma}  \label{lem:compacts_in_triangulated}
  Let~$\GEN$ be a countable set of compact objects of a triangulated
  category~$\Tri$ that has countable direct sums.  Let
  $A\in\gen{\GEN}$.  Then~$A$ belongs to $\gen{\GEN}_\fin$ if and only
  if~$A$ is compact.
\end{lemma}

\begin{proof}
  Objects of $\gen{\GEN}_\fin$ are compact because the compact objects
  in a triangulated category form a thick triangulated subcategory.
  Up to the fact that we only have countable direct sums, the converse
  is a result of Amnon Neeman \cite{Neeman:Localising}*{Lemma 2.2}.
  We have to check that his proof only uses countable direct sums
  under our cardinality hypotheses.  This is routine, so that we omit
  the verification.  We mention that the critical points of the
  argument are explained in greater detail
  in~\cite{Neeman:Grothendieck_duality}.
\end{proof}

\begin{lemma}  \label{lem:CIZ_compactly_generated}
  Let~$G$ be a second countable locally compact group.  If $H\subseteq
  G$ is smooth, then $C_0(G/H)$ is compact.  Up to conjugacy there are
  only countably many smooth compact subgroups.
\end{lemma}

\begin{proof}
  Let $N\subseteq H\subseteq U\subseteq G$ be such that~$U$ is open
  and almost connected and $U/N$ is a Lie group.  Let $K\subseteq U$
  be a maximal compact subgroup containing~$H$.
  \cite{MeyerNest}*{Corollary 3.2} implies $\KK^G_0(C_0(G/H),B) \cong
  \KK^K_0(C_0(U/H),\Res_G^K B)$ because $C_0(G/H)\cong \Ind_U^G
  C_0(U/H)\cong \Ind_K^G C_0(K/H)$.  We want to obtain a
  $K$\nbd{}equivariant isomorphism $C_0(U/H)\cong C_0(U/K)\times
  C(K/H)$.  Then the compactness of $C_0(U/H)$ follows immediately
  from the Poincaré duality isomorphism (see \cite{Kasparov})
  \begin{displaymath}
    \KK^K_0(C_0(U/H),B)
    \cong \KK^K_0(\C,C_0(U/K) \otimes C_\tau(K/H)\otimes B).
  \end{displaymath}

  We prove $U/H\cong U/K\times K/H$ following Herbert Abels
  (\cite{Abels:Parallel}).  The Lie algebra~$\mathfrak{k}$ of $K/N$
  acts on the Lie algebra of~$U/N$ by conjugation.  We split the
  latter into invariant subspaces $\mathfrak{k}\oplus\bigoplus_{i=1}^n
  \mathfrak{p}_i$.  Abels observes that the map
  $$
  \prod \mathfrak{p}_i \times K \to U,
  \qquad
  \bigl((x_i),k\bigr) \mapsto \prod \exp(x_i) \cdot k
  $$
  is a $K$\nbd{}equivariant diffeomorphism.  This yields the desired
  $K$\nbd{}equivariant diffeomorphism $U/H\cong \prod
  \mathfrak{p}_i\times K/H \cong U/K\times K/H$.
\end{proof}

\begin{prop}  \label{pro:ADir_finite}
  Let~$G$ be a second countable locally compact group and let~$\ADir$
  be the source of the Dirac morphism.  Then
  $\ADir\in\gen{\CI_1}_\fin$ if and only if the functor $B\mapsto
  \RKK^G_0(\EG;\C,B)$ commutes with countable direct sums and only
  produces countable groups for separable~$B$.  If~$G$ is almost
  totally disconnected, this is also equivalent to
  $\ADir\in\gen{\CI_0}_\fin$.

  A sufficient condition for this is the existence of a
  $G$\nbd{}compact model for $\EG$.
\end{prop}

\begin{proof}
  Lemma~\ref{lem:CIZ_compactly_generated} shows that
  Lemma~\ref{lem:compacts_in_triangulated} applies to $\CI_1$, and to
  $\CI_0$ if~$G$ is almost totally disconnected.  We have
  $\ADir\in\gen{\CI_1}$ by \cite{MeyerNest}*{Proposition 9.2}, and
  $\ADir\in\gen{\CI_0}$ if~$G$ is almost totally disconnected.
  Theorem~\ref{the:derived_category} and Lemma~\ref{lem:RKKG_as_Ktheory}
  yield
  \begin{equation}  \label{eq:KK_ADir}
    \KK^G_0(\ADir,B)
    \cong \RKK^G_0(\EG;\C,B)
    \cong \K_0(C_0(\abs\EG,B)\cross G).
  \end{equation}
  Thus $\ADir$ is a compact object of $\KK^G$ if and only if
  $\RKK^G_0(\EG;\C,B)$ has the properties required in the statement of
  the theorem.  Moreover, if we can find a $G$\nbd{}compact model for
  $\EG$, then $C_0(\abs\EG)$ is a $C^*$-algebra.  In this case, the isomorphism 
\eqref{eq:KK_ADir} implies immediately that $\ADir$ is compact by
 well-known properties of $\K$-theory.
\end{proof}

We can now prove Theorem~\ref{the:dual_Dirac_finite}.  We remark that
the argument only uses the potentially weaker hypothesis
$\ADir\in\gen{\CI_1}_\fin$.  However, we know no example of a group 
satisfying the latter condition but without a
 $G$\nbd{}compact model for $\EG$.

\begin{proof}[Proof of Theorem~\ref{the:dual_Dirac_finite}]
  By Theorem~\ref{the:equivariant_versus_coarse}, the assumption on
  $\mu_{\abs{G},H}^*$ implies that
  \begin{equation}  \label{eq:equivariant_coassembly_B}
    p_\EG^*\colon \KK^G(\C,B)\to\RKK^G(\EG;\C,B)
  \end{equation}
  is an isomorphism for all $B\in\CI_1$.  The class of objects~$B$ for
  which~\eqref{eq:equivariant_coassembly_B} is an isomorphism is a
  thick triangulated subcategory of $\KK^G$.  Therefore, it contains
  $\gen{\CI_1}_\fin$.  By Proposition~\ref{pro:ADir_finite}, it
  contains~$\ADir$.  Thus a dual Dirac morphism exists by
  Lemma~\ref{lem:dual_Dirac_isomorphism}.  The converse assertion is
   Theorem~\ref{the:dual_Dirac_gives_isomorphism}.
\end{proof}

\subsection{Discrete groups with finite dimensional classifying space}
\label{sec:finite_dim_EG}

We now pass to the case when~$G$ has a finite \emph{dimensional} model
for $\EG$.  We also assume~$G$ to be \emph{discrete}.  Therefore,
$\EG$ can be realised by a simplicial complex.  We assume that
\begin{enumerate}[(1)]
\item we can choose $\EG$ finite dimensional;

\item there are only finitely many conjugacy classes of finite
  subgroups in~$G$.

\end{enumerate}
Of course, the second condition is trivially satisfied for torsion
free groups.  It is known that there exist groups with finite
dimensional $\EG$ that violate (2).

As in Section~\ref{sec:finite_EG}, what we really need is a condition
on the domain of the Dirac morphism~$\ADir$.  Let $\CI_2$ be the set
of $C^*$-algebras of the form 
$C_0(\N\times G/H)$ as $H$ ranges over the
 finite subgroups of $ G$.  We let
$\gen{\CI_2}_\fin$ be the thick triangulated subcategory generated by
$\CI_2$.  Since $C_0(G/H)$ is a retract of $C_0(\N\times G/H)$, this
category contains $\gen{\CI_0}_\fin$.  The same argument as in the
proof of Theorem~\ref{the:dual_Dirac_finite} yields the following
lemma:

\begin{lemma}  \label{lem:dual_Dirac_finite_dim}
  Suppose that $\ADir\in\gen{\CI_2}_\fin$.  Then~$G$ has a dual Dirac
  morphism if and only if the $H$\nbd{}equivariant coarse
  co\brd{}assembly map with coefficients $C_0(\N)$ is an isomorphism
  for all finite subgroups $H\subseteq G$.
\end{lemma}

\begin{thm}  \label{bigtheorem2}
  Let~$G$ be a discrete group that satisfies the conditions (1)
  and (2) above.  Then $\ADir\in\gen{\CI_2}_\fin$.  Hence~$G$ has a
  dual Dirac morphism if and only if the $H$\nbd{}equivariant coarse
  co\brd{}assembly map with coefficients $C_0(\N)$
  $$
  \mu_{\abs{G},H,C_0(\N)}^*\colon
  \K_{*+1}(\bvr(\abs{G},C_0(\N))\cross H) \to \KX^*_H(\abs{G},C_0(\N))
  $$
  is an isomorphism for all finite subgroups $H\subseteq G$.
\end{thm}

\begin{proof}
  We use constructions of Gennadi Kasparov and Georges Skandalis
  in~\cite{KasparovSkandalis}.  Choose a finite dimensional simplicial
  model for $\EG$ as above.  By a barycentric subdivision, we can
  arrange that the action of~$G$ is ``type preserving''.  Hence~$G$
  acts on the $C^*$-algebra $\mathcal{A}_\EG$ defined
  in~\cite{KasparovSkandalis}.  As we observed in
  Corollary~\ref{cor:Dirac_EG_simplicial}, $\mathcal{A}_\EG$ is a
  model for~$\ADir$.  Thus we have to prove that
  $\mathcal{A}_\EG\in\gen{\CI_2}_\fin$.  In different notation, this
  is already shown in~\cite{KasparovSkandalis}.  We sketch the
  argument.  The skeletal filtration of $\EG$ gives rise to a
  filtration of $\mathcal{A}_\EG$ by ideals
  $$
  0 = \mathcal{A}^{(-1)}_\EG \subseteq
  \mathcal{A}^{(0)}_\EG \subseteq
  \mathcal{A}^{(1)}_\EG \subseteq
  \mathcal{A}^{(2)}_\EG \subseteq \dots \subseteq
  \mathcal{A}^{(n)}_\EG =
  \mathcal{A}_\EG,
  $$
  where~$n$ is the dimension of $\EG$.  The resulting extensions
  $$
  0 \to
  \mathcal{A}^{(k-1)}_\EG \to
  \mathcal{A}^{(k)}_\EG \to
  \mathcal{A}^{(k)}_\EG/\mathcal{A}^{(k-1)}_\EG \to
  0
  $$
  have $G$\nbd{}equivariant, completely positive, contractive sections
  because all occurring $G$-$C^*$\nbd{}algebras are proper and
  nuclear.  Hence these extensions are admissible. By induction
  on~$k$, it follows that $\mathcal{A}^{(k)}_\EG$ belongs to the
  triangulated subcategory of $\KK^G$ generated by the subquotients
  $\mathcal{A}^{(j)}_\EG/\mathcal{A}^{(j-1)}_\EG$ for all~$j$.  Thus
  it remains to prove that these subquotients belong to
  $\gen{\CI_2}_\fin$.  They are $\KK^G$\nbd{}equivalent to
  $C_0(\EG^{(j)})$, where $\EG^{(j)}$ denotes the set of
  $j$\nbd{}cells of~$\EG$, viewed as a discrete $G$\nbd{}space.  Thus
  $\EG^{(j)}$ is a disjoint union of homogeneous spaces $G/H$ for
  finite subgroups $H\subseteq G$.  By assumption, there are at most
  finitely many non-isomorphic proper homogeneous spaces $G/H$.  Hence
  we can write $\EG^{(j)}$ as a \emph{finite} disjoint union of spaces
  of the form $G/H\times I$, where~$I$ is some countable set.  This
  implies $C_0(\EG^{(j)})\in\gen{\CI_2}_\fin$ as desired.
\end{proof}

\begin{cor}
  If~$G$ is a torsion free discrete group with finite dimensional
  classifying space $\EG$, then~$G$ has a dual Dirac morphism if and
  only if
  $$
  \mu_{\abs{G},C_0(\N)}^*\colon
  \K_{*+1}(\bvr(\abs{G},C_0(\N))) \to \KX^*(\abs{G},C_0(\N))
  $$
  is an isomorphism.  In particular, the existence of a dual Dirac
  morphism for~$G$ is a coarse property, that is, it only depends on
  the coarse space~$\abs{G}$.
\end{cor}

\section{Geometric K-theory}
\label{sec:Kgeo}

In this section we define the closely related notions of
\emph{boundary class} and \emph{regularising} class in
$\RKK^G(\EG;\C,\C)$ and observe that such classes give rise to
homotopy invariant higher signatures.  The boundary classes arise from
the equivariant topological $\K$\nbd{}theory of the stable Higson
corona.  The regularising classes arise from the Dirac dual Dirac
method.  For totally disconnected groups with $G$\nbd{}compact $\EG$,
boundary classes and regularising classes coincide.  We will also use
these ideas to construct dual Dirac morphisms and boundary classes
from certain geometric situations.

\subsection{Boundary classes and regularising classes}
\label{sec:boundary_classes}

\begin{prop}  \label{pro:nat_trafo_Kreg}
  Let $\Dirac\in\KK^G(\ADir,\C)$ be a Dirac morphism for~$G$.  Then
  the following diagram commutes 
  $$
  \xymatrix{
    \Ktop_{*+1}\bigl(G,\bvr_G(\abs{G})\bigr)
    \ar@/^3em/[rr]^-{\mu^*_{\abs{G},G}}
    \ar[d]^{\Psi_*} \ar[r]^-{\partial_{\abs{G},\C,\ADir}} &
    \K_*(C_0(\abs\EG,\ADir)\cross G) \ar[d]_{\cong}
    \ar[r]^{\Dirac_*}_{\cong} &
    \K_*(C_0(\abs\EG)\cross G) \ar[d]_{\cong} \\
    \KK_*^G(\C,\ADir) \ar[r]^-{p_\EG^*} \ar[dr]_-{\Dirac_*} &
    \RKK_*^G(\EG;\C,\ADir) \ar[r]^{\Dirac_*}_{\cong} &
    \RKK_*^G(\EG;\C,\C) \\
    & \KK^G(\C,\C) \ar[ur]_-{p_\EG^*},
  }
  $$
  where~$\Psi_*$ is as in Corollary \ref{cor:nat_trafo} and 
  $\mu^*_{\abs{G},G}$ is the $G$\nbd{}equivariant coarse co-assembly
  map for~$\abs{G}$. 
\end{prop}

\begin{proof}
  The left square diagram is the special case of
  Corollary~\ref{cor:nat_trafo} with $B=\C$ and $P=\ADir$, together
  with our identification $\Ktop_*(G,A) = \K_*\bigl((A\otimes
  \ADir)\cross G\bigr)$.  The rest of the diagram evidently commutes,
  given our definition of the coarse co-assembly map.  The vertical
  isomorphisms in the right square are contained in
  Lemma~\ref{lem:RKKG_as_Ktheory}.  The other isomorphisms follow from
  the invertibility of $p_\EG^*(\Dirac)$, see
  Section~\ref{sec:resolutions_Dirac_dual_Dirac}.
\end{proof}

The above diagram points to a particular subgroup of
$\RKK^G_*(\EG;\C,\C)$, namely the range of the composition
$\Dirac_*\circ p_\EG^*= p_\EG^*\circ\Dirac_*\colon \KK^G_*(\C,\ADir)
\to \RKK^G_*(\EG;\C,\C)$.  It contains the classes that come from
$\Ktop_{*+1}\bigl( G, \bvr(\abs{G})\bigr)$ via $\mu_{\abs{G},G}$.
Accordingly we make the following definition.

\begin{defn}  \label{def:boundary_regularising}
  Let $\alpha\in\RKK^G_*(\EG;\C,\C)$.  We call~$\alpha$ a
  \emph{boundary class} if it belongs to the range of the map
  $\Ktop_{*+1}\bigl(G,\bvr_G(\abs{G})\bigr)\to\RKK^G_*(\EG;\C,\C)$ in
  Proposition~\ref{pro:nat_trafo_Kreg}.  We call~$\alpha$
  \emph{regularising} if it belongs to the range of
  the map $\KK^G_*(\C,\ADir) \to 
  \RKK^G_*(\EG ; \C , \C)$ of the same Proposition. 
\end{defn}

It follows immediately from Proposition~\ref{pro:nat_trafo_Kreg} that
boundary classes are regularising and that regularising classes belong
to the range of $p_\EG^*\colon \KK^G_*(\C,\C)\to\RKK^G_*(\EG;\C,\C)$
and hence yield homotopy invariant higher signatures. The group
$\Ktop_{*+1}\bigl(G,\bvr_G(\abs{G})\bigr)$ should be thought of as a
sort of \emph{geometric $\K$-theory group} for~$G$. Given its 
importance, we state this as a proposition.

\begin{prop} \label{pro:homotopyinvarianceofboundaryclasses}
Regularising classes, whence also boundary classes, 
are in the range of $p_\EG^*$ and thus yield homotopy invariant
higher signatures.
\end{prop}

We also have the following.

\begin{prop}  \label{pro:ideals}
  The boundary classes and the regularising classes form graded ideals
  in $\RKK^G_*(\EG;\C,\C)$.  
\end{prop}

Thus, if $\alpha,\beta\in\RKK^G_*(\EG;\C,\C)$ and~$\alpha$ is
regularising, then $\alpha\beta$ is also regularising and hence
homotopy invariant.  This is the reason for calling these classes
regularising.

\begin{proof}
  We return to our discussion of the diagrams \eqref{eq:KK_diagram}
  and~\eqref{eq:RKK_diagram}.  We have already observed that the maps
  $$
  \KK^G(\ADir,\ADir) \overset{p_\EG^*}\longrightarrow
  \RKK^G(\EG;\ADir,\ADir) \overset{\tau_\ADir}\longleftarrow
  \RKK^G(\EG;\C,\C)
  $$
  are isomorphisms.  Since they are graded ring homomorphisms,
  $\KK^G(\ADir,\ADir)$ and $\RKK^G(\EG;\C,\C)$ are isomorphic as
  graded rings.  It is well-known that the graded ring
  $\RKK^G(\EG;\C,\C)$ is graded
  commutative.  Hence there is no difference between left and right
  modules over it.
  
  If~$F$ is a functor defined on $\KK^G$, then $\KK^G(\ADir,\ADir)$
  automatically acts on $F(A\otimes\ADir)$.  In particular, it acts on
  the graded Abelian groups $\Ktop(G)$, $\KK^G(\C,\ADir)$, and
  $\RKK^G(\EG;\C,\ADir)$.  Consider the left square in the diagram in
  Proposition~\ref{pro:nat_trafo_Kreg}.  Except for~$\Psi_*$, the
  other three maps are evidently module homomorphisms.  This implies
  that the images of $\Ktop\bigl(G,\bvr(G)\bigr)$ and
  $\KK^G(\C,\ADir)$ in $\RKK^G(\EG;\C,\C)$ are ideals.
\end{proof}

We are compelled to conclude on purely formal grounds
 that regularising classes and
boundary classes are rather special even among classes in the range of
$p_\EG^*$, and in particular amongst those classes yielding
homotopy invariant higher signatures.  
The unit is always among the latter.  On the other hand, if it is
regularising or a boundary class, then \emph{every} class is such
because these form ideals.  In this case $G$ actually has a 
dual Dirac morphism by Lemma~\ref{lem:dual_Dirac_isomorphism}.

We also note that if we use the traditional definition of $\Ktop(G)$,
then the action of $\RKK^G_*(\EG;\C,\C)$ on topological
$\K$\nbd{}theory discussed above may be described rather easily.  Thus
the $\RKK^G_*(\EG;\C,\C)$\brd{}module structure of $\Ktop_*\bigl(G;
\bvr(\abs{G})\bigr)$ is not special to our setup.  However, the
construction of the isomorphism~$\Psi_*$ identifying regularising and
boundary classes is difficult with the traditional definition.

If $X\in\gen{\CI}$, in particular if~$X$ is proper, then $1_X \otimes \Dirac
$ is invertible in $\KK^G(X\otimes \ADir , X)$.  Given $A,B\in\KK^G$,
$X\in\gen{\CI}$, and $\alpha\in\KK^G(A,X)$, $\beta\in\KK^G(X,B)$, we
define $\beta\bullet\alpha\in\KK^G(A,B\otimes\ADir)$ as the
composition $(1_\ADir \otimes \beta) \circ (1_X\otimes \Dirac )^{-1} \circ \alpha
\in \KK^G(A, \ADir \otimes B)$. 
%$$
%A \overset{\alpha}\to X \overset{\Dirac^{-1}_*}\longrightarrow
%X\otimes\ADir \overset{\beta_*}\longrightarrow B\otimes\ADir.
%$$
Naturality of exterior products implies that
$\Dirac_*(\beta\bullet\alpha)=\beta\circ\alpha$. Hence 
$p_\EG^*( \beta \circ \alpha)  = p_\EG^* \circ \Dirac_*(  \beta 
\bullet \alpha )  = \Dirac_* \circ p_\EG^* ( \beta \bullet \alpha)$ is 
regularising. Conversely, by the definition of regularising, 
any regularising class has the form $\Dirac_*(\alpha)$ for 
some $\alpha \in \RKK^G(\C , \ADir)$.

% and
%$\Dirac^*(\beta\bullet\alpha)=\tau_\ADir(\beta\circ\alpha)$.  If
%$A=B=\C$, it follows that $\beta\circ\alpha$ is regularising.
%Conversely, any regularising element of $\RKK^G(\EG;\C,\C)$ has the
%form $\Dirac\circ\alpha$ for some $\alpha\in\KK^G(\C,\ADir)$.

To summarize, an element of $\RKK^G(\EG;\C,\C)$ is regularising if and
only if it can be factorised through an object of $\gen{\CI}$.  For
many groups, $\ADir$ is a proper $G$\nbd{}$C^*$\brd{}algebra, so that
elements of~$I$ even factor through a proper
$G$\nbd{}$C^*$\brd{}algebra.

The following is a central result of this paper. It may be regarded as
a refinement of Theorem \ref{the:dual_Dirac_finite}. It implies in 
particular that for groups $G$ with $G$-compact $\EG$, 
regularising classes and boundary classes 
are the same

\begin{thm}  \label{the:Kreg}
  Let~$G$ be an almost totally disconnected group with
  $G$\nbd{}compact $\EG$. Then for every $B \in \KK^G$, the map 
  $$
  \Psi^B_*\colon \Ktop_{*+1}\bigl(G,\bvr(\abs{G},B)\bigr) \to
  \KK_*^G(\C,B\otimes\ADir)
  $$
  of Corollary~\ref{cor:nat_trafo} is an isomorphism for which 
  the diagram
  $$
  \xymatrix{
    \Ktop_{*+1}\bigl(G,\bvr(\abs{G},B)\bigr)
    \ar[r]^-{\mu_{\abs{G},G,B}^*}
    \ar[d]_{\cong}^{\Psi^B_*} &
    \KX^*_G(\abs{G},B) \ar[d]_{\cong} \\
    \KK_*^G(\C,B\otimes\ADir) \ar[r]^-{p_\EG^*} &
    \RKK_*^G(\EG;\C,B\otimes\ADir)
  }
  $$
  commutes.  In particular, $\Ktop_{*+1}\bigl(G,\bvr(\abs{G})\bigr)$
  is naturally isomorphic to $\KK^G_*(\C,\ADir)$. 

  Moreover, $\KK_*^G(\C,B\otimes\ADir)\cong\KK_*^G(\C,B)$ if~$B$ is a
  proper $G$\nbd{}$C^*$\brd{}algebra or, more generally, if
  $B\in\gen{\CI}$.
\end{thm}

\begin{proof}
  Corollary~\ref{cor:nat_trafo} with $P=\ADir$ yields the desired
  diagram.  The map $\Psi^{B,P}_*$ is an isomorphism for $P\in\CI_0$.
  The class of~$P$ for which this is the case is thick and
  triangulated.  Hence it contains~$\ADir$ by
  Proposition~\ref{pro:ADir_finite}. The isomorphism statement follows. 
  Proper $G$\nbd{}$C^*$\brd{}algebras belong to
  $\gen{\CI}$ by \cite{MeyerNest}*{Corollary 7.3}.  We have
  $B\otimes\ADir\cong B$ if and only if $B\in\gen{\CI}$ by
 \cite{MeyerNest}*{Theorem 4.7}.
\end{proof}

We summarize our main result as follows. 

\begin{cor}   \label{cor:detect_regularising}
  Let~$G$ be an almost totally disconnected group with
  $G$\nbd{}compact $\EG$ and let $a\in\RKK^G_*(\EG;\C,\C)$.  Then the
  following are equivalent:
  \begin{enumerate}
  \item $a$ can be factorised through an object of $\gen{\CI}$;
  \item $a$ is regularising;
  \item $a$ is a boundary class.
  \end{enumerate}
\end{cor}

\subsection{Dual Dirac morphisms from compactifications}
\label{sec:compactify}

We are going to construct dual Dirac morphisms from contractible
admissible compactifications, thereby strengthening a result of Nigel Higson
(\cite{Hig2}).  We assume throughout that there is a $G$\nbd{}compact
model for $\EG$.  Hence there is no difference between $C_0(\EG)$ and
$C_0(\abs\EG)$.

Recall that a metrisable compactification $Z\supseteq\abs\EG$ is
called \emph{admissible} if all (scalar valued) continuous functions
on~$Z$ have vanishing variation.  If this is the case for scalar
valued functions, it automatically holds for operator valued functions
because $C(Z,D)\cong C(Z)\otimes D$.  An \emph{equivariant
  compactification} of $\abs{\EG}$ is a compactification together with
a $G$\nbd{}action that extends the given action on $\abs\EG$.

\begin{thm}  \label{the:dual_Dirac_compactify}
  Let~$G$ be a locally compact group with a $G$\nbd{}compact model for
  $\EG$ and let $\abs{\EG}\subseteq Z$ be an admissible equivariant
  compactification.  If~$Z$ is $H$\nbd{}equivariantly contractible for
  all compact subgroups $H\subseteq G$, then~$G$ has a dual Dirac
  morphism.
\end{thm}

\begin{proof}
  Since~$Z$ is admissible, we have an embedding
  $\bbbarred_G(Z)\subseteq\vvr_G(\abs{\EG})$.  Let $\partial Z\defeq
  Z\setminus \EG$ be the boundary of the compactification.
  Identifying
  $$
  \bbbarred(\partial Z)\cong\bbbarred(Z)/C_0(\EG,\Comp_G),
  $$
  we obtain a morphism of extensions
  $$
  \xymatrix{
    0 \ar[r] & C_0(\abs\EG,\Comp_G) \ar[r] \ar@{=}[d] &
    \vvr_G(\abs\EG) \ar[r] \ar[d]^{\subseteq} &
    \bvr_G(\abs\EG) \ar[r] \ar[d]^{\subseteq} & 0 \\
    0 \ar[r] & C_0(\abs\EG,\Comp_G) \ar[r] &
    \bbbarred_G(Z) \ar[r] &
    \bbbarred_G(\partial Z) \ar[r] & 0.
  }
  $$
  Let~$H$ be a compact subgroup.  Since $Z$ is compact, $\bbbar_G(Z)$
  is identical with $C(Z,\Comp)$.  In particular, since~$Z$ is
  $H$\nbd{}equivariant contractible by hypothesis, $\bbbar_G(Z)$ is
  $H$\nbd{}equivariantly homotopy equivalent to $\C$.  Hence
  $\bbbarred_G(Z)$ has vanishing $H$\nbd{}equivariant
  $\K$\nbd{}theory.  This implies
  $\Ktop_*\bigl(G,\bbbarred_G(Z)\bigr)= 0$ by~\cite{Ect}, so that the
  connecting map
  $$
  \Ktop_{*+1}\bigl(G,\bbbarred_G(\partial Z)\bigr)
  \to \Ktop_*\bigl(G,C_0(\EG)\bigr)
  \cong \K_*(C_0(\EG)\cross G)
  $$
  is an isomorphism.  This in turn implies that the connecting map
  $$
  \Ktop_{*+1}\bigl(G,\bvr_G(\abs\EG)\bigr)
  \to \Ktop_*\bigl(G,C_0(\EG)\bigr)
  $$
  is surjective.  Thus we can lift $1\in \RKK^G_0(\EG;\C,\C) \cong
  \K_0(C_0(\EG)\cross G)$ to
  $$
  \alpha\in \Ktop_1\bigl(G,\bvr_G(\abs\EG)\bigr)
  \cong \Ktop_1\bigl(G,\bvr_G(\abs{G})\bigr).
  $$
  Then $\Psi_*(\alpha)\in \KK^G_0(\C,\ADir)$ is the desired dual
  Dirac morphism.
\end{proof}

Theorem~\ref{the:dual_Dirac_compactify} applies to the Gromov boundary
for a hyperbolic group.  We conclude that a dual Dirac morphism must
come from the topological $\K$\nbd{}theory of the Gromov boundary.  In
fact, the above argument shows that all boundary classes for a
hyperbolic group come from the Gromov boundary.

\section{Lipschitz and proper Lipschitz K-theory classes}
\label{sec:Lipschitz}

The idea of~\cite{CGM2} is to prove homotopy invariance for a
particular higher signature by showing that it arises from a specific
geometric construction.  This motivated us to formulate the notion of
boundary classes and prove Theorem~\ref{the:coarse_gives_classes},
which is a direct analogue of some constructions in~\cite{CGM2}.  In
this section, we show how Lipschitz cohomology classes and boundary
classes are related.  In the case of proper Lipschitz classes, the
relationship is exact: every proper Lipschitz class is a boundary
class. The proof is a consequence of 
Corollary~\ref{cor:detect_regularising}: proper Lipschitz classes
are regularising. 

 We also use the stable Higson corona construction to approach
non proper Lipschitz classes, simplifying the geometric part of the
proof of homotopy invariance of Gelfand-Fuchs cohomology classes, as well 
as to describe a more general geometric situation in which the essential idea 
of~\cite{CGM2} can be made to work.

\subsection{Construction of higher signatures}
\label{sec:constructionofhighersignatures}

Let $X$ be a locally compact $G$-space, and suppose that $G$ acts
by translations with respect to some coarse structure on 
 $X$. From Theorem~\ref{the:coarse_gives_classes}
we have a map 
$$ \K_{*+1}\bigl(\bvr(X)\bigr)
\overset{\psi_*^{G,X}}\longrightarrow \KK^G_*\bigl(\C,C_0(X)\bigr).
$$
 Hence
given any nonzero class in the $\K$-theory of the stable
Higson corona of $X$, we may push it forward and obtain a 
class $\alpha \in \KK^G_N(\C , C_0 (X))$. This class in turn induces a 
map $\K_*\bigl( \Cmax (G)\bigr) \to \K_{*+N}\bigl( C_0 (X) \cross G)\bigr)$. 
If
 $b\colon \K_{*+N}(C_0(X)\cross G)\to\R$ is a linear functional, then
 we may then construct a higher signature for $G$ (see the Introduction 
for a brief discussion of higher signatures), by 
 the composition
\begin{equation}\label{constructinofoahighersig}
\Ktop_*(G) \to \K_*(\Cmax G)
\overset{\alpha_*}\to \K_{*+N}(C_0(X) \cross G) \overset{b}\to \R,
\end{equation}
where the first map is the analytic assembly map.  This higher
signature is homotopy invariant by construction since it factorises
through the analytic assembly map.

Such a `geometric' higher signature therefore has two components; 
the construction of a nonzero class in the $\K$-theory of the Higson 
corona of a $G$-space $X$, and the construction of a linear functional 
$b$ as above. The former is obviously related to coarse geometry; the
latter is, however, not. If $X$ is a \emph{proper} $G$-space, then 
topology provides us with many classes in
 $\KK^G_*(C_0 (X), \C)$, whence linear functionals $b$ above, 
and the nontrivial part of the above procedure is 
finding nonzero classes in the $\K$-theory of the stable Higson
corona of $X$. It is possible to find such classes if the 
coarse co-assembly map for $X$ is an isomorphism, for then 
the $\K$-theory of the stable Higson corona is the same 
as the coarse $\K$-theory of $X$.

If $X$ is proper, then the higher signature constructed 
above corresponds to a \emph{regularising class} in $\RKK^G_*(\EG ; \C ,\C)$. 
Indeed, asuming that our linear functional $b$ comes from 
a class 
 $\beta \in \KK^G(C_0 (X), \C)$, we can form the class
 $ \beta\bullet\alpha  \in\KK^G(\C, \ADir)$. We may then 
apply the isomorphism $\Dirac_* \circ p_\EG^* : \KK^G(\C , \ADir)
\stackrel{\cong}{\to} \RKK^G_*(\EG ; \C , \C)$ to obtain a class
 in $ \RKK^G_*(\EG ; \C , \C)$. The latter class is by construction 
\emph{regularising}. The higher signature corresponding to 
it (see Introduction) is the same as that constructed by the 
recipe~\eqref{constructinofoahighersig}. If $G$ has a $G$-compact
model for $\EG$, we conclude that the higher signatures 
constructed by~\eqref{constructinofoahighersig} where $X$ is a 
proper $G$-space correspond in fact to \emph{boundary classes}.

In this paper we discuss three natural examples where the above situation
arises for an action of $G$ on a locally compact space $X$: 
when $X$ has the coarse structure coming from an admissible 
compactification; when $X$ admits commuting actions of $G$ and another
group $H$, whence admits a coarse structure inherited from the action of $H$; 
and when $X$ admits a proper map to Euclidean space 
satisfying a certain displacement condition, whence admits 
a pulled-back coarse 
structure. 

 The first 
situation we have already implicitly discussed whilst
proving Theorem~\ref{the:dual_Dirac_compactify}. The second 
will appear in the context of Gelfand-Fuchs classes below. The 
third 
amounts to a reformulation of the ideas of Connes, Gromov and 
Moscovici in~\cite{CGM2}, and we begin with it.

\subsection{Pulled-back coarse structures and Lipschitz classes}
\label{sec:pulled_back}

%More precisely, suppose as above that $X$ is a \emph{proper} $G$-space.  In this case, for
%most purposes we may assume that the linear functional~$b$ comes from
%an element $\beta\in \KK^G_*(C_0(X),\C)\cong \KK(C_0(X)\cross G,\C)$.
%We may then form the class $p_\EG^*(\beta \circ [\alpha])\in
%\RKK^G_*(\EG;\C,\C)$, and it is easy to check that the resulting
%higher signature
%$$
%\Ktop_*(G) \to \Z \subset \R,
%\qquad a\mapsto \mathrm{ind}_G(a\cdot p_\EG^*(\beta\circ[\alpha]),
%$$
%agrees with the one defined above.  Since~$X$ is assumed proper, the
%class $p_\EG^*(\beta\circ[\alpha])$ is always regularising by the
%construction of $\beta\bullet[\alpha]$.  Using
%Corollary~\ref{cor:detect_regularising}, we get:%%%

%\begin{cor}
%  If~$G$ is a discrete group with a $G$\nbd{}compact model for $\EG$,
%  then every proper Lipschitz $\K$\nbd{}theory class in
%  $\RKK^G_*(\EG;\C,\C)$ is a boundary class.
%\end{cor}

Let~$X$ be a $G$\nbd{}space, let~$Y$ be a coarse space and let
$\alpha\colon X\to Y$ be a proper continuous map.  We pull back the
coarse structure on~$Y$ to a coarse structure on~$X$, so that
$E\subseteq X\times X$ is an entourage if and only if
$\alpha_*(E)\subseteq Y\times Y$ is one.  Since~$\alpha$ is proper and
continuous, this coarse structure is compatible with the topology
on~$X$.  The group~$G$ acts by translations with respect to this coarse 
structure if and only if~$\alpha$
satisfies the following \emph{displacement condition}: for any compact
subset $K\subseteq G$,
$$
\bigl\{\bigl(\alpha(gx),\alpha(x)\bigr) \bigm|
x\in X,\ g\in K\bigr\} \subseteq Y\times Y
$$
is an entourage.  The map~$\alpha$ becomes a coarse map. Hence we obtain canonical maps
$$
\K_{*+1}\bigl(\bvr(Y)\bigr)
\overset{\alpha^*}\longrightarrow \K_{*+1}\bigl(\bvr(X)\bigr)
\overset{\psi_*^{X,G}}\longrightarrow \KK^G_*\bigl(\C,C_0(X)\bigr).
$$
where $\psi_*^{G,X}$ is as in Theorem~\ref{the:coarse_gives_classes}.

In particular, 
we can push forward any class in the group $\K_{*+1}\bigl(\bvr(Y)\bigr) $ 
to a class $\alpha \in\KK^G_*\bigl(\C,C_0(X)\bigr)$
and proceed as in Section~\ref{sec:constructionofhighersignatures}.

The constructions of~\cite{CGM2}*{Section I.10} involve the 
specific choice of
 $Y=\R^N$ with the coarse structure from
the Euclidean metric.  The coarse co\nbd{}assembly map is an
isomorphism for~$\R^N$ because $\R^N$ is scalable.  Moreover, $\R^N$
is uniformly contractible and has bounded geometry.  Hence the map
$j\colon \R^N\to\Rips(\R^N)$ is a coarse homotopy equivalence. We
have, therefore canonical
isomorphisms
$$
\K_{*+1}\bigl(\bvr(\R^N)\bigr)
\cong \K_*\bigl(C_0(\R^N)\bigr)
\cong \K_{*+N}(\C)
$$ and thus the $\K$-theory of $\bvr(\R^N)$ is infinite cyclic.  We
denote the generator by $[\partial\R^N]
\in\K_{1-N}\bigl(\bvr(\R^N)\bigr)$.  Of course, this is nothing but
the usual dual Dirac morphism for the locally compact group $\R^N$.
Thus we obtain a class
$$
[\alpha]\defeq \psi\bigl(\alpha^*[\partial\R^N]\bigr)
\in\KK^G_{-N}\bigl(\C,C_0(X)\bigr)
$$
for any map $\alpha\colon X\to\R^N$ that satisfies the displacement
condition above.

There is a slightly more general setup, also contained in~\cite{CGM2},
where we replace a map to~$\R^N$ by a section of an
$N$\nbd{}dimensional vector bundle, as follows.  Let~$P$ be a
$G$\nbd{}space and let $\pi\colon X\to P$ be a $G$\nbd{}equivariant
$\Spin(N)$\nbd{}principal bundle.  That is, the action of~$G$ on~$X$
commutes with the action of $\Spin(N)$ and~$\pi$ is
$G$\nbd{}equivariant.  Let $T\defeq X\times_{\Spin(N)} \R^N$ be the
associated vector bundle over~$P$.  It carries a $G$\nbd{}invariant
Euclidean metric and spin structure.  If $\alpha\colon P\to T$ is a
section, then we can define a map $\alpha'\colon X\to\R^N$ by sending
$x\in X$ to the coordinates of $\alpha\pi(x)$ in the orthogonal frame
described by~$x$.  This map is $\Spin(N)$\brd{}equivariant with
respect to the standard action of $\Spin(N)$ on $\R^N$.  Conversely,
any $\Spin(N)$\brd{}equivariant map $\alpha'\colon X\to\R^N$ arises in
this fashion.  Since $\Spin(N)$ is compact, the map~$\alpha'$ is
proper if and only if $p\mapsto \norm{\alpha(p)}$ is a proper function
on~$P$.

As above, we can use a $\Spin(N)$\nbd{}equivariant proper continuous
map $\alpha'\colon X\to\R^N$ to pull back the coarse structure of
$\R^N$ to~$X$.  Then $\Spin(N)$ acts by isometries.  The group~$G$
acts by translations if and only if the section $\alpha\colon P\to T$
associated to~$\alpha'$ satisfies the displacement condition that
$$
\sup \{ \norm{g\alpha(g^{-1} x)-\alpha(x)} | x\in X,\ g\in K\}
$$
be bounded for all compact subsets $K\subseteq G$.  Suppose
that~$\alpha$ satisfies this.  Then we are in the situation of
Theorem~\ref{the:coarse_gives_classes} with $H=\Spin(N)$.  Since~$H$
acts freely, $C_0(X)\cross\Spin(N)$ is Morita-Rieffel equivalent to
$C_0(P)$, and this Morita-Rieffel equivalence is $G$\nbd{}equivariant.
We obtain canonical maps
\begin{multline*}
  \K^{\Spin(N)}_{*+1}\bigl(\bvr_{\Spin(N)}(\R^N)\bigr)
  \overset{(\alpha')^*}\longrightarrow
  \K^{\Spin(N)}_{*+1}\bigl(\bvr_{\Spin(N)}(X)\bigr)
  \\ \overset{\psi}\longrightarrow
  \KK^G_*\bigl(\C,C_0(X)\cross \Spin(N)\bigr)
  \cong \KK^G_*\bigl(\C,C_0(P)\bigr).
\end{multline*}

The space~$\R^N$ is $\Spin(N)$\brd{}equivariantly scalable, so that
the $\Spin(N)$\nbd{}equivariant coarse co-assembly map for~$\R^N$ is
an isomorphism.  Moreover, the action of $\Spin(N)$ on $\R^N$ is spin
by definition.  Hence
$$
\K^{\Spin(N)}_{*+1}\bigl(\bvr_{\Spin(N)}(\R^N)\bigr)
\cong \K^{\Spin(N)}_*\bigl(C_0(\R^N)\bigr)
\cong \K^{\Spin(N)}_{*+N}(\C).
$$
The usual dual Dirac morphism $[\partial
\R^N]\in\K^{\Spin(N)}_{1-N}\bigl(\bvr_{\Spin(N)}(\R^N)\bigr)$ for
$\R^N$ is the image of the trivial representation of $\Spin(N)$ in
$\K^{\Spin(N)}_0(\C)$.  As above, we obtain a class $[\alpha]\defeq
\psi\bigl(\alpha^*[\partial\R^N]\bigr) \in\KK^G_{-N}\bigl(\C,C_0(X)\bigr)$ for
any proper section $\alpha\colon P\to T$ satisfying the
displacement condition.

The classes constructed in the above manner
 are called \emph{Lipschitz classes}; if the space $X$ is 
proper, they are called proper Lipschitz classes. 
In view of our discussion
 in Section~\ref{sec:constructionofhighersignatures}, we have:

\begin{cor}
  If~$G$ is a discrete group with a $G$\nbd{}compact model for $\EG$,
  then every proper Lipschitz $\K$\nbd{}theory class in
  $\RKK^G_*(\EG;\C,\C)$ is regularising, whence a boundary class.
\end{cor}

\subsection{Coarse structures on jet bundles}
\label{sec:Gelfand-Fuchs}

We recall the setup for the proof of homotopy invariance of
Gelfand-Fuchs cohomology classes in~\cite{CGM2}.  Let~$M$ be an
oriented compact manifold and let $\Diff^+(M)$ be the infinite
dimensional Lie group of orientation preserving diffeomorphisms
on~$M$.  Let the locally compact group~$G$ act on~$M$ by orientation
preserving diffeomorphisms, that is, by a continuous group
homomorphism $G\to\Diff^+(M)$.  We are interested in the classes in
the cohomology of~$G$ that we obtain by pulling back the Gelfand-Fuchs
cohomology for~$M$ ( the latter being part of the group cohomology of
$\Diff^+(M)$; see \cite{CGM2}).

Let $\pi^k\colon \Jet^k_+(M)\to M$ be the \emph{oriented $k$\nbd{}jet
  bundle} over~$M$.  That is, a point in $\Jet^k_+(M)$ is the $k$th
order Taylor series at~$0$ of a germ of an orientation preserving
diffeomorphism of a neighbourhood of $0\in\R^n$ into~$M$.  Germs of
diffeomorphisms of neighbourhoods of $0\in\R^n$ form a connected Lie
group~$H$.  Its Lie algebra~$\mathfrak{h}$ is the space of polynomial
maps $p\colon \R^n\to\R^n$ of order~$k$ with $p(0)=0$, with an
appropriate Lie algebra structure.  The maximal compact subgroup
$K\subseteq H$ is isomorphic to $\mathrm{SO}(n)$, acting by isometries
on~$\R^n$.  It acts on~$\mathfrak{h}$ by conjugation.

The bundle $\Jet^k_+(M)$ is an $H$\nbd{}principal bundle over~$M$.
Since the action of~$H$ on $\Jet^k_+(M)$ is natural, it commutes with
the action of~$G$.  We let~$H$ act on the right and~$G$ on the left.
Define $X_k\defeq \Jet^k_+(M)/K$.  This is the bundle space of a
fibration over~$M$ with fibres $H/K$.

The complex that computes Gelfand-Fuchs cohomology can be represented
in a canonical way as a complex of $\Diff^+(M)$\brd{}invariant
differential forms on~$X_k$.  These differential forms generate cyclic
cocycles on $C^\infty_{\mathrm{c}}(X_k)\cross_{\mathrm{alg}} G$.
Using the formalism of $n$\nbd{}traces, these cyclic cocycles can be
extended to linear functionals $\K_*(C_0(X_k)\rcross G)\to\C$.  Once
we have an appropriate element $\alpha\in\KK^G(\C,C_0(X_k))$, we can
pull back these classes to linear functionals on $\K_*(\Cred G)$ and
thus prove the homotopy invariance of Gelfand-Fuchs cohomology
classes.  We can construct~$\alpha$ in the following fashion.

We equip $\Jet^k_+(M)$ with the unique coarse structure for which~$H$
acts isometrically defined in Theorem~\ref{the:proper_act_coarse}.
The compactness of $\Jet^k_+(M)/H\cong M$ implies easily that~$G$ acts
by translations on $\Jet^k_+(M)$.  Moreover, any orbit map $\abs{H}\to
\Jet^k_+(M)$ is a coarse equivalence.  We have a Morita-Rieffel
equivalence $C_0(X_k)\sim C_0(\Jet^k_+ M)\cross K$ because~$K$ acts
freely on $\Jet^k_+(M)$.  Thus we want to look at the map
$$
\K_{*+1}(\bvr_K(\Jet^k_+ M)\cross K) \overset{\psi}\to
\KK^G_*(\C,C_0(\Jet^k_+ M)\cross K) \cong
\KK^G_*\bigl(\C,C_0(X_k)\bigr)
$$
produced by Theorem~\ref{the:coarse_gives_classes}.  Since~$H$ is
almost connected, it has a dual Dirac morphism by~\cite{Kasparov}.
Moreover, $H/K$ is a model for $\EG$ by~\cite{Abels:Parallel}.
Together with Theorem~\ref{the:dual_Dirac_gives_isomorphism} this
implies
$$
\K_{*+1}(\bvr_K(\Jet^k_+ M)\cross K)
\cong \K_{*+1}(\bvr_K(\abs{H})\cross K)
\cong \KX^*_K(\abs{H})
\cong \K^*_K(H/K).
$$
Let $\mathfrak{h}$ and~$\mathfrak{k}$ be the Lie algebras of $H$
and~$K$.  There is a $K$\nbd{}equivariant homeomorphism
$\mathfrak{h}/\mathfrak{k}\cong H/K$, where~$K$ acts on
$\mathfrak{h}/\mathfrak{k}$ by conjugation.  Now we need to know
whether there is a $K$\nbd{}equivariant spin structure on
$\mathfrak{h}/\mathfrak{k}$.  One can check that this is the case
if~$n$ is even or if $k\equiv 0,1 \bmod 4$.  Since we can choose~$k$
as large as we like, we can always assume that this is the case.
Hence we may identify $\K^*_K(H/K)$ with the representation ring
of~$K$ in degree~$-N$, where $N=\dim \mathfrak{h}/\mathfrak{k}$.  The
trivial representation of~$K$ yields a canonical element in
$\KK^G_{-N}\bigl(\C,C_0(X_k)\bigr)$.

The above construction is essentially the same as in~\cite{CGM2}.  We
get some additional trouble with spin structures because we want to
use only the case of ``fixed target'' in the notation of~\cite{CGM2}.
Our framework allows us to use the existence of a dual Dirac morphism
for~$H$.  In contrast, this fact is reproved in different notation
in~\cite{CGM2}.

\end{document}